\documentclass[reqno]{amsart}

\usepackage[numbers, sort]{natbib}
\usepackage[letterpaper, total={6.5 in, 9 in}, centering]{geometry}  
\usepackage{amsmath}
\usepackage{amssymb}
\usepackage[overload]{empheq}
\usepackage{amsfonts}
\usepackage{enumerate}
\usepackage{amsthm}
\usepackage{setspace}
\usepackage{multirow}
\usepackage{multicol}
\usepackage{mathdots}
\usepackage{mathtools}
\usepackage{graphicx}
\usepackage{float}
\usepackage{mathrsfs}
\usepackage{cases}
\usepackage[all]{xy}
\usepackage{subcaption}
\usepackage{chemarrow}
\usepackage{extarrows}
\usepackage{mathrsfs}
\usepackage{bm}
\usepackage{units}
\usepackage{pgfplots}
\usepackage{pgfplotstable}
\usepgfplotslibrary{external}
\pgfplotsset{compat=1.9,
	log x ticks with fixed point/.style={
		xticklabel={
			\pgfkeys{/pgf/fpu=true}
			\pgfmathparse{exp(\tick)}%
			\pgfmathprintnumber[fixed relative, precision=3]{\pgfmathresult}
			\pgfkeys{/pgf/fpu=false}
		}
	},
	log y ticks with fixed point/.style={
		yticklabel={
			\pgfkeys{/pgf/fpu=true}
			\pgfmathparse{exp(\tick)}%
			\pgfmathprintnumber[fixed relative, precision=3]{\pgfmathresult}
			\pgfkeys{/pgf/fpu=false}
		}
	}
}

\usepackage{scalerel,stackengine}
\stackMath
\newcommand\reallywidehat[1]{%
\savestack{\tmpbox}{\stretchto{%
  \scaleto{%
    \scalerel*[\widthof{\ensuremath{#1}}]{\kern-.6pt\bigwedge\kern-.6pt}%
    {\rule[-\textheight/2]{1ex}{\textheight}}
  }{\textheight}%
}{0.5ex}}%
\stackon[1pt]{#1}{\tmpbox}%
}

\newtheorem{lemma}{Lemma}

\usepackage{natbib}

\DeclareMathOperator*{\ext}{ext}
\DeclareMathOperator*{\dexp}{d \exp}
\DeclareMathOperator*{\d2exp}{d^2 \exp}
\DeclareMathOperator*{\ad}{ad}
\DeclareMathOperator*{\tr}{tr}
\DeclareMathOperator*{\pol}{pol}
\DeclareMathOperator*{\sym}{Sym}
\DeclareMathOperator*{\sk}{Skew}
\DeclareMathOperator*{\asym}{Asym}
\DeclareMathOperator*{\lyap}{Lyap}

\author{Xuefeng Shen}
\author{Khoa Tran}
\author{Melvin Leok}
\title {High-order symplectic Lie group methods on $SO(n)$ using the polar decomposition}

\begin{document}

\begin{abstract}
	A variational integrator of arbitrarily high-order on the special orthogonal group $SO(n)$ is constructed using the polar decomposition and the constrained Galerkin method. It has the advantage of avoiding the second-order derivative of the exponential map that arises in traditional Lie group variational methods. In addition, a reduced Lie--Poisson integrator is constructed and the resulting algorithms can naturally be implemented by fixed-point iteration. The proposed methods are validated by numerical simulations on $SO(3)$ which demonstrate that they are comparable to variational Runge--Kutta--Munthe-Kaas methods in terms of computational efficiency. However, the methods we have proposed preserve the Lie group structure much more accurately and and exhibit better near energy preservation.
\end{abstract}


\maketitle

\section{Introduction}
	
	\subsection{Overview}
	Given a configuration manifold $Q$, variational integrators provide a useful method of deriving symplectic integrators for Lagrangian mechanics on the tangent bundle $TQ$ in terms of the Lagrangian $L$, or for Hamiltonian mechanics on the cotangent bundle $T^*Q$ in terms of the Hamiltonian $H$. It involves discretizing Hamilton's principle or Hamilton's phase space principle rather than the Euler--Lagrange or Hamilton's equations. Discrete Lagrangian variational mechanics is described in terms of a discrete Lagrangian $L_d(q_0,q_1)$, which is an approximation of the exact discrete Lagrangian,
	\begin{equation}\label{ldexact}
		L_d^{\text{exact}}(q_0,q_1) 
		= \ext_{\substack{q\in C^2([0,h],Q) \\ q(0)=q_0, \text{ } q(h)=q_1}} \int_0^h L(q(t), \dot q(t)) \, dt.
	\end{equation}
	The discrete Hamilton's principle states that the discrete action sum is stationary for variations that vanish at the endpoints. This yields the  discrete Euler--Lagrange equations,
	\[
		D_2 L_d(q_{k-1},q_k) + D_1 L_d(q_k,q_{k+1}) = 0,
	\]
	where $D_i$ denotes the partial derivative with respect to the $i$-th argument. This defines an update map on $Q \times Q$, where $(q_{k-1},q_k) \mapsto (q_k,q_{k+1})$. The update map can equivalently be described in terms of the discrete Legendre transforms,
	\begin{equation}
	\label{vi}
		p_k =  -D_1 L_d(q_k, q_{k+1}), \qquad
		p_{k+1} =  D_2 L_d(q_k, q_{k+1}),
	\end{equation}
	which defines an update map $(q_k,p_k)\mapsto (q_{k+1},p_{k+1})$ on $T^*Q$ that automatically preserves the canonical symplectic structure on $T^*Q$.
	
	The order of the variational integrator depends on how accurately $L_d(q_0,q_1)$ approximates $L_d^{\text{exact}}(q_0,q_1)$. To derive a high-order discrete Lagrangian, a typical approach is the Galerkin method~\cite{MaWe2001}. This involves considering the definition of the exact discrete Lagrangian, replacing $C^2([0,h],Q)$ with a finite-dimensional function space, and approximating the integral with a numerical quadrature formula. When the configuration manifold $Q$ is a linear space and the polynomials of degree less than or equal to $s$ are chosen, the classical symplectic partitioned Runge--Kutta methods are recovered. Subsequently, \citet{Leok2012} introduced the shooting-based discrete Lagrangian, which allows one to construct a symplectic integrator from an arbitrary one-step method. 

	When the configuration manifold $Q$ is a Lie group $G$, the construction of the discrete Lagrangian is more complicated than the case of linear space. \citet{Leok2005} proposed parametrizing curves on the Lie group using the exponential map, namely a curve $g(t)$ connecting $g_0$ and $g_1$ that is represented by
	\[
		g(t) = g_0\cdot \exp(\epsilon(t)),
	\]
	where $\epsilon(t)\in \mathfrak{g}$ is a curve on the Lie algebra of $G$ with fixed endpoints $\epsilon(0) = 0$ and $\epsilon(h)= \log(g_0^{-1} g_1)$. This allows one to replace variations in $g(t)$ by variations in $\epsilon(t)$ on the Lie algebra $\mathfrak{g}$, which is a linear space. This yields the following expression for the exact discrete Lagrangian,
	\begin{equation}
	\label{ldgroup}
		L_d^{\text{exact}}(g_0,g_1) 
		= \ext_{\substack{\epsilon\in C^2([0,h],\mathfrak{g}) \\ \epsilon(0) = 0, \text{ } \epsilon(h) = \log(g_0^{-1}g_1)}} 
		\int_0^h L(g_0\cdot \exp(\epsilon(t)), g_0\cdot \dexp \nolimits_{\epsilon(t)} (\dot{\epsilon}(t))) \, dt,
	\end{equation}
	where $\dexp_{\epsilon}(\dot{\epsilon}) = \exp(\epsilon)\cdot \frac{1 - \text{e}^{-\text{ad}_{\epsilon}}} {\text{ad}_{\epsilon}} (\dot{\epsilon})$ is the tangent lift of the exponential map. If $\epsilon(t)$ is restricted to a finite-dimensional function space and the integral is replaced with a quadrature rule, we obtain the Galerkin Lie group variational integrators. The error analysis and implementation details for such methods can be found in \cite{Hall2014,Geir2014}. The above construction can be naturally extended to any retraction~\cite{abmase2008} on $G$, which is a diffeomorphism from a neighborhood of $0\in \mathfrak{g}$ to neighborhood of $e\in G$ that satisfies a rigidity condition. The main disadvantage of Galerkin Lie group variational integrators is the term $\dexp$ in \eqref{ldgroup}. This implies that the resulting discrete Euler--Lagrange equations involve $\d2exp$, which cannot be calculated exactly in general and requires the truncation of a series expansion.

	\subsection{Contributions}
	In this paper, we focus on the Lie group $SO(n)$ as our configuration space. By using the fact that every invertible square matrix can be uniquely decomposed into the product of an orthogonal matrix and a symmetric positive-definite matrix via the polar decomposition, we will circumvent the disadvantage discussed previously: Instead of parametrizing curves on $SO(n)$ by the exponential map or a retraction, $SO(n)$ is embedded naturally in the space $GL_{+}(n) = \{A\in \mathbb{R}^{n\times n}\mid \text{det}(A)>0\}$, an open subset of $\mathbb{R}^{n\times n}$. Given fixed endpoints $g_0$ and $g_1$, we will construct interpolating polynomials in $GL_{+}(n)$ while ensuring that the internal points remain on $SO(n)$ by using the polar decomposition. Furthermore, we do not extend the Lagrangian $L$ to $GL_{+}(n)$ but instead project the trajectory onto $SO(n)$ in the same way. The variational integrator in Lagrangian form is derived following the usual variational approach for the constrained Galerkin method and the Hamiltonian form is derived using the discrete Legendre transforms.
	
	For a system with rotational symmetry, we obtain a simpler integrator using Lie--Poisson reduction on the Hamiltonian side. Namely, if $L$ is $SO(n)$-invariant, the constructed discrete Lagrangian is also $SO(n)$-invariant and we can construct a reduced symplectic Lie--Poisson integrator. Lastly, we consider the dipole on a stick problem from \citep{Geir2014}, conduct numerical experiments using our method, and compare these to the  variational Runge--Kutta--Munthe-Kaas methods (VRKMK) from the same reference. 

\section{Background}
\label{sect:Background}	
	\subsection{Notation}
	\label{sect:BackgroundNotations}
	Recall that the Lie algebra of $SO(n)$ is the set $\sk (n) = \{\Omega \in \mathbb{R}^{n\times n} \mid \Omega^{T} = -\Omega\}$, with the matrix commutator as the Lie bracket. Here, the inner products on $\mathbb{R}^{n\times n}$ and $\sk (n)$ are introduced, and we identify these spaces with their dual spaces using the Riesz representations. For any $ A,B\in \mathbb{R}^{n\times n}$, the inner product is given by
	\[
		\tr \left(A B^T\right) = \sum\limits_{i,j=1}^n a_{ij} b_{ij}.
	\]
	For any $\Omega, \tilde{\Omega} \in \sk (n)$, the inner product is defined by
	\[
		\langle \Omega, \tilde{\Omega} \rangle = \sum\limits_{i<j}\Omega_{ij} \tilde{\Omega}_{ij} = \frac{1}{2} \tr (\Omega \tilde{\Omega}^T).
	\]
	In addition, consider the operator $\asym \colon \mathbb{R}^{n \times n} \to \sk(n)$, defined by $\asym(A) = A - A^T$.
	The following properties can be easily verified:
	\begin{enumerate}[(a)]
		\item For any $A,B\in \mathbb{R}^{n\times n}$, $\tr (AB^T)= \tr (A^{T}B)$.
		\item For any $A,B,P,Q\in \mathbb{R}^{n\times n}$, $\tr (A(PBQ)^T) = \tr((P^TAQ^T)B^T)$.
		\item For any $\Omega \in \sk (n), A\in \mathbb{R}^{n\times n}$, $\langle\Omega, \asym (A)\rangle = \sum\limits_{i<j}\Omega_{ij}(A_{ij}-A_{ji}) = \tr (\Omega A^T)$.
	\end{enumerate}
	In particular, we note that (c) gives a relationship between the two inner products. Lastly, given the choice of inner products, Riesz representation allows us to identify $(\mathbb{R}^{n\times n})^*$ with $\mathbb{R}^{n\times n}$ and $\sk (n)^*$ with $\sk (n)$.
	\subsection{Polar Decomposition}
	\label{sect:BackgroundPD}
	We introduce the \textit{polar decomposition} and the construction of the retraction on $SO(n)$ described in \cite{celledoni2002class}. Given $A\in GL(n)$, it decomposes uniquely as
	\begin{equation*}
		A \longrightarrow UP, \qquad U \in O(n), \text{ } P \in {\sym}^+(n),
	\end{equation*}
	 where ${\sym}^+(n)$ is the set of $n \times n$ symmetric positive-definite matrices. This is the polar decomposition of $A$, and we denote it as a coproduct mapping by $\pol \colon GL(n) \to O(n) \times {\sym}^+(n)$. The map of interest is the projection $\mathbb{P} = \pi_1 \circ \pol \colon GL(n) \to O(n)$ defined by
	 \begin{equation*}
	 	\mathbb{P}(A) = U,
	 \end{equation*}
	 where $\pi_1$ is the projection onto the first coordinate. In particular, when $A\in GL_{+}(n)$, we have $U\in SO(n)$. A fast and efficient algorithm for calculating the projection of the polar decomposition is by Newton iteration,
	\begin{equation}
	\label{newtoniterate}
		U_{k+1} = \frac{1}{2}\left(U_k + U_k^{-T}\right), \text{ } U_0 = A.
	\end{equation}
	
	Now, this projection can be used to construct a retraction on $SO(n)$ from its Lie algebra $\sk(n)$ relative to the identity element $I$,
	\begin{equation*}
		\mathbb{P}(I + \Omega) = U,
	\end{equation*}
	where $\Omega \in \sk (n)$. This provides a diffeomorphism between a neighborhood of $0\in \sk (n)$ and a neighborhood of $I \in SO(n)$. To calculate its inverse, suppose that $I + \Omega = UP$ and take the transpose on both sides to obtain $I - \Omega = PU^{T}$, which implies that $U^{T}(I+\Omega) = (I-\Omega)U$. Thus, we have that
	\begin{equation}
	\label{pinverse}
   		U^{T}\Omega + \Omega U + U^{T}-U=0.
	\end{equation}
This is a\textit{ Lyapunov equation}, and it is well-known that matrix equations of the form $AX+XB+C=0$ have a unique solution if and only if for any $\lambda \in \sigma(A), \mu \in \sigma(B)$, $\lambda+\mu \neq 0$. For $U$ in the neighborhood of identity, its eigenvalues lie in the open right-half plane, which ensures that a unique solution to \eqref{pinverse} exists. In principle, this Lyapunov equation can be solved using classical algorithms \cite{bartels1972, golub1979}. For convenience, we denote the solution to the Lyapunov equation as
	\[
		X = \text{Lyap}(A,B,C).
	\]
	Next we introduce the tangent map and its adjoint for $\mathbb{P}$, which are essential for the derivation of the variational integrator.

	\subsubsection{The Tangent Map}
	\label{sect:BackgroundPDTM}
	Consider the polar decomposition $A(t) = U(t)P(t)$ and differentiate both sides to yield $\dot{A} = \dot{U}P+U\dot{P}$. By left-trivialization on $SO(n)$, we can write $\dot{U} = U\Omega$, where $\Omega \in \sk (n)$. Rearranging gives $\dot{P} = U^T\dot{A} - \Omega P$, and since   $\dot{P} \in {\sym}^+(n)$, we get that $U^{T}\dot{A}-\Omega P = \dot{A}^{T}U+P\Omega$. Consequently, we may write it in the form of a Lyapunov equation, 
	\begin{equation}
	\label{dP}
    	P\Omega + \Omega P + \dot{A}^{T}U -  U^{T}\dot{A} = 0. 
	\end{equation}
	We see that the tangent map of the polar decomposition $d\mathbb{P}_A: \mathbb{R}^{n\times n}\to \sk (n)$ is given by $d\mathbb{P}_A(\dot{A}) = \Omega$, where we solve the Lyapunov equation \eqref{dP} to obtain $\Omega$. 
	
	\subsubsection{The Adjoint of the Tangent Map}
	\label{sect:BackgroundPDAdjointTM}
	The adjoint of $d\mathbb{P}_A$ can be defined as
	\[
		\tr (d \mathbb{P}_A^*(\Omega) B^T) = \langle\Omega, d \mathbb{P}_A(B)\rangle,\qquad \forall \Omega\in \sk (n), B \in \mathbb{R}^{n\times n}.
	\]
	Recall that $d \mathbb{P}_A(B)$ involves solving the Lyapunov equation \eqref{dP}. We aim to compute $d \mathbb{P}_A^*$, so we define the following two maps,
	\begin{alignat*}{2}
    	\phi	& \colon \sk (n) \longrightarrow \sk (n),			&	\qquad \Omega	& \longmapsto \Omega P + P\Omega, \\
    	\psi	& \colon \mathbb{R}^{n\times n} \longrightarrow \sk (n),	&	\qquad B 		& \longmapsto U^{T}B - B^{T}U,
	\end{alignat*}
	where $A = UP$ is fixed. Therefore, $d \mathbb{P}_A$ can be viewed as composition of $\psi$ and $\phi^{-1}$, 
	\[
		d \mathbb{P}_A = (\phi^{-1}\circ \psi) \colon \mathbb{R}^{n\times n} \longrightarrow \sk (n),
	\]
and $d \mathbb{P}_A^*(\Omega) = (\phi^{-1}\circ \psi)^*(\Omega) = \psi^*\circ (\phi^*)^{-1}(\Omega)$. We shall derive the expressions for $\phi^*$ and $\psi^*$ by considering the Riesz representations for our domains and codomains. For the adjoint of $\phi$, let  $ \Omega, \tilde{\Omega}\in \sk (n)$, then
	\begin{align*}
    	\langle\phi^*(\Omega), \tilde{\Omega} \rangle
    	&= \langle \Omega, \phi(\tilde{\Omega})\rangle = \langle \Omega, \tilde{\Omega}P + P\tilde{\Omega} \rangle\\
		&= \langle \Omega, \asym (\tilde{\Omega}P)\rangle = \tr (\Omega(\tilde{\Omega}P)^T) 	\\
		&= \tr ((\Omega P)\tilde{\Omega}^T) \\
		&= \langle\Omega P+P\Omega, \tilde{\Omega} \rangle.
	\end{align*}
Thus, $\phi^* = \phi$, and $\phi$ is Hermitian. For $\psi^*$, let $\Omega\in \sk (n)$ and $B \in \mathbb{R}^{n\times n}$, then
	\begin{align*}
    	\tr (\psi^*(\Omega) B^T) 
		&= \langle \Omega, \psi(B)\rangle = \langle\Omega, U^{T}B - B^{T}U\rangle \\
		&= \langle \Omega, \asym (U^{T}B) \rangle = \tr (\Omega (U^TB)^T) \\
		&= \tr ((U\Omega)B^T).
\end{align*}
	Therefore, $\psi^*(\Omega) = U\Omega$, and we obtain
	\[
		d \mathbb{P}_A^*(\Omega) = (\psi^*\circ (\phi^*)^{-1})(\Omega) = (\psi^*\circ \phi^{-1})(\Omega) = U \lyap (P,\Omega^{T}),
	\]
where $\text{Lyap}(P,\Omega^{T}) = \text{Lyap}(P,P,\Omega^{T})$. Finally, we state a lemma that will be used later:
	\begin{lemma}
	\label{null}
		$\mathbb{P}(I+S) = I$ if and only if $S\in \sym (n)$ and $\lambda > -1$ for all $\lambda \in \sigma(S)$.
	\end{lemma}
	

\section{Lagrangian Variational Integrators on the Rotation Group $SO(n)$}
\label{sect:VI}
	Let the configuration manifold be the rotation group $G = SO(n)$, and $L \colon G \times \mathfrak{g} \to \mathbb{R}$ be a left-trivialized Lagrangian. We shall construct a discrete Lagrangian following the approach of constrained Galerkin methods on $GL_+(n)$ (see Appendix \ref{append_constrainedGalerkin}). Denote the internal points by $\{U_i\}_{i=1}^s \subset G$ and the left-trivialized internal tangent vectors by $\{\Omega_i\}_{i=1}^s \subset \mathfrak{g}$. Fixing the endpoints $g_0$ and $g_1$, we have 
	\[
		L_d(g_0,g_1) = h \sum_{i=1}^s b_i L(U_i,\Omega_i),
	\]
subject to the constraint
	\begin{equation}
	\label{restrict}
		g_1 = \mathbb{P}\left(g_0 + h\sum\nolimits_{i=1}^s b_i U_i \Omega_i \right),
	\end{equation}
where the internal points $U_i$ are represented by
	\begin{equation}
	\label{internal}
		U_i = \mathbb{P}\left(g_0 + h\sum\nolimits_{j=1}^s a_{ij} U_j \Omega_j \right).
	\end{equation}
	The expressions inside the parentheses in \eqref{restrict} and \eqref{internal} correspond to a Runge--Kutta method in the embedding space. But, since these points may not lie on the Lie group $G$, they are projected to $G$ using the polar decomposition.

	Observe that \eqref{restrict} is equivalent to the condition that $\mathbb{P}\left(g_1^{T}(g_0 + h\sum_{i=1}^s b_i U_i \Omega_i) \right) = I$. Suppose that $h$ is small, and $g_0$ and $g_1$ are close enough to each other, then $g_1^{T}\left(g_0 + h\sum_{i=1}^s b_i U_i \Omega_i \right)$ is in the neighborhood of $I$. By Lemma \ref{null}, \eqref{restrict} holds if and only if $g_1^{T} \left(g_0 + h \sum_{i=1}^s  b_i U_i \Omega_i \right)\in \text{Sym}(n)$, and so it is equivalent to
	\[
		\asym \left(g_1^{T} \left(g_0 + h\sum\nolimits_{i=1}^s b_i U_i \Omega_i \right)\right) = 0.
	\]
	Now, we can construct a discrete Lagrangian with the constraint using a Lagrange multiplier  $\Lambda \in \sk (n)$,
	\[
		\tilde{F} (g_0,g_1, \{\Omega_i\}_{i=1}^s, \Lambda) 
		= h \sum_{i=1}^s b_i L(U_i,\Omega_i) + \left\langle \Lambda, \asym \left( g_1^{T} \left(g_0 + h \sum\nolimits_{i=1}^s b_i U_i \Omega_i \right) \right) \right\rangle.
	\]
	The corresponding variational integrator is given by
	\begin{subequations}
	\begin{align}[left = \empheqlbrace\,]
 		0 &= \frac{\partial \tilde{F}}{\partial \Omega_i}, \qquad i=1,2\dots s,
 		\label{omega}\\
 		0 &= \frac{\partial \tilde{F}}{\partial \Lambda}, 
 		\label{lambda} \\
		U_i &= \mathbb{P}\left(g_0 + h \sum\nolimits_{j=1}^s a_{ij} U_j \Omega_j\right),\\
		-p_0 &= \frac{\partial \tilde{F}}{\partial g_0}, 
		\label{g0}\\
		p_1 &= \frac{\partial \tilde{F}}{\partial g_1}. 
		\label{g1}
	\end{align}
	\end{subequations}
	
	
	We shall compute the above equations more explicitly. It is easy to see that \eqref{lambda} is equivalent to the constraint \eqref{restrict}. Let us turn to \eqref{omega}, where the main difficulty is the implicit dependence of $\{U_i\}_{i=1}^s$ on $\{\Omega_i\}_{i=1}^s$ that involves solving a nonlinear system \eqref{internal}. Suppose that $k\in \{1, 2, \ldots, s\}$ is fixed, and we vary $\Omega_k$ such that $\Omega_k \to \Omega_k(t)$ with $\Omega_k(0) = \Omega_k$ and $\dot{\Omega}_k(0) = \delta\Omega_k$, while $\{\Omega_i\}_{i\neq k}$ remain fixed. Then,
	\[
		U_i(t) = \mathbb{P} \left( g_0 + h\sum\nolimits_{j\neq k} a_{ij} U_j(t) \Omega_j + h a_{ik} U_k(t) \Omega_k(t) \right).
	\]
	Differentiating both sides and letting $\dot{U}_i = U_iX_{ik}$, we have that
	\begin{equation}
	\label{temp1}
		X_{ik} = d \mathbb{P}_{A_i} \left( h \sum\nolimits_{j=1}^s a_{ij} U_j X_{jk} \Omega_j + h a_{ik} U_k \delta\Omega_k \right),
	\end{equation}
	where $A_i = g_0 + h \sum\nolimits_{j=1}^s a_{ij} U_j \Omega_j$. Then, \eqref{temp1} can be rewritten as 
	\begin{equation}
	\label{temp2}
		X_{ik} - d \mathbb{P}_{A_i} \left(h \sum\nolimits_{j=1}^s a_{ij} U_j X_{jk} \Omega_j\right) = h d \mathbb{P}_{A_i}(a_{ik} U_k \delta\Omega_k).
	\end{equation}
	In order to represent $\{X_{ik}\}_{i=1}^s$ in terms of $\delta\Omega_k$, we define three maps,
	\begin{alignat*}{3}
		\tilde{\psi}_k	&: \sk (n) \longrightarrow \sk (n)^s, 	&	\quad	\delta \Omega_k & \longmapsto \{d \mathbb{P}_{A_i}(a_{ik} U_k \delta \Omega_k) \}_{i=1}^s, \\
		\tilde{\phi} 	&: \sk (n)^s \longrightarrow \sk (n)^s,	&	\quad 	\{X_{ik}\}_{i=1}^s & \longmapsto \left\{ X_{ik} - d \mathbb{P}_{A_i} \left( h \sum\nolimits_{j=1}^s a_{ij} U_j X_{jk} \Omega_j \right) \right\}_{i=1}^s,\\
		\tilde{\pi}_i 	&: \sk (n)^s \longrightarrow \sk (n),	&	\quad 	\{\Omega_i\}_{i=1}^s & \longmapsto \Omega_i.
	\end{alignat*}
	Then, 
	\begin{equation}
	\label{Xik}
		X_{ik} = (\tilde{\pi}_i \circ \tilde{\phi}^{-1} \circ h\tilde{\psi}_k)(\delta\Omega_k) = h(\tilde{\pi}_i \circ \tilde{\phi}^{-1} \circ \tilde{\psi}_k)(\delta\Omega_k).
	\end{equation}
	Now, we compute $\frac{\partial \tilde{F}}{\partial \Omega_k}$ by evaluating $\left.\frac{d}{d t}\right|_{t=0} \tilde{F}(g_0, g_1, \ldots, \Omega_k(t), \ldots, \Lambda)$ and expressing $\frac{\partial L}{\partial U}: G\times \mathfrak{g}\to \mathfrak{g}^*$ as a left-trivialized cotangent vector. Since this is a straightforward calculation of the variation, we present the result here for equation \eqref{omega} (see section \ref{append_dFdOmega} for details),
	\begin{equation}
	\label{nomega}
		0 = h \sum_{i=1}^s b_i (\tilde{\pi}_i \circ \tilde{\phi}^{-1}\circ \tilde{\psi}_k)^* \left( \frac{\partial L}{\partial U} (U_i, \Omega_i) + \asym (U_i^{T} g_1 \Lambda\Omega_i^{T}) \right) +  b_k \left( \frac{\partial L}{\partial \Omega}(U_k,\Omega_k) + \asym (U_k^{T} g_1 \Lambda) \right),
	\end{equation}
	for any $k = 1, 2, \ldots, s$. 
	
	Recall that $\tilde{\pi}_i \circ \tilde{\phi}^{-1}\circ \tilde{\psi}_k \colon \sk (n) \to \sk (n)$ and its dual is given by
	\[
		(\tilde{\pi}_i \circ \tilde{\phi}^{-1}\circ \tilde{\psi}_k)^*
		= 
		\tilde{\psi}_k^* \circ (\tilde{\phi}^*)^{-1} \circ \tilde{\pi}_i^*.
	\]
	Therefore, let us derive the explicit expressions for the adjoints of our three proposed maps, so we may write $(\tilde{\pi}_i \circ \tilde{\phi}^{-1}\circ \tilde{\psi}_k)^*$ explicitly. The adjoint of $\tilde{\pi}_i$ is easy to compute, and for any $S \in \sk (n)$,
	 \begin{equation}
	 \label{pistar}
    	\tilde{\pi}_i^*(S) = (0, \ldots, S, \ldots 0),
	\end{equation}
	where $S$ is in the $i$-th position. For the adjoint of $\tilde{\phi}$, we consider the identity
	\[
		\langle \tilde{\phi}^*(S_1, \ldots, S_s), (\tilde{S}_1, \ldots, \tilde{S}_s) \rangle
		= \langle (S_1, \ldots, S_s), \tilde{\phi}(\tilde{S}_1, \ldots, \tilde{S}_s) \rangle,
	\]
	for any $(S_1, \ldots, S_s), (\tilde{S}_1, \ldots, \tilde{S}_s) \in \sk (n)^s$. Using the properties of the inner products again, we obtain the explicit expression (see section \ref{append_phistar} for details),
	\begin{equation}
	\label{phistar}
		\tilde{\phi}^*(S_1, \ldots, S_s) 
		= \left\{ S_j - \asym \left(h U_j^T \sum\nolimits_{i=1}^s a_{ij} d \mathbb{P}_{A_i}^*(S_i) \Omega_j^T  \right)\right\}_{j=1}^{s}.
	\end{equation}
	 Similarly, we consider the same identity and techniques to obtain the explicit expression for $\tilde{\psi}_k^*$ (see section \ref{append_psistar} for details),
	\begin{equation}
	\label{psistar}
		\tilde{\psi}_k^*(S_1, \ldots, S_s) = \asym\left( U_k^T \sum\nolimits_{i=1}^s a_{ik} d \mathbb{P}_{A_i}^*(S_i) \right).
	\end{equation}
	Combining \eqref{pistar}, \eqref{phistar}, and \eqref{psistar},  $(\tilde{\pi}_i \circ \tilde{\phi}^{-1}\circ \tilde{\psi}_k)^*(S)$ for $S \in \sk (n)$ can be computed as follows,
	\begin{subequations}
	\label{ff}
	\begin{align}[left = \empheqlbrace\,]
    	S_j - \asym \left( h U_j^{T} \sum\nolimits_{l=1}^s a_{lj} d \mathbb{P}_{A_l}^*(S_l) \Omega_j^{T} \right) &= S \cdot \delta_{ij},	\qquad j=1,2\dots s,	\label{f1}	\\
    	(\tilde{\pi}_i \circ \tilde{\phi}^{-1}\circ \tilde{\psi}_k)^*(S) &= \asym \left( U_k^{T} \sum\nolimits_{l=1}^s a_{lk} d\mathbb{P}_{A_l}^*(S_l)\right). \label{f2}
	\end{align}	
	\end{subequations}
	We can first calculate $\{S_l\}_{l=1}^s$ from \eqref{f1} by using fixed-point iteration, and then substitute the result into \eqref{f2} to obtain $(\tilde{\pi}_i\circ \tilde{\phi}^{-1}\circ \tilde{\psi}_k)^*(S)$. This combined with \eqref{nomega} gives an explicit formula for \eqref{omega}.
	
	We now derive an explicit formula for \eqref{g0}. Notice that $U_i$ depends on $g_0$ by the nonlinear system in \eqref{internal}, so we can use the method of variations again. If we vary $g_0 \rightarrow g_0(t)$ such that $g_0(0) = g_0$ and $\dot{g}_0(0) = g_0 \delta g_0$, we obtain
	\[
		U_i(t) = \mathbb{P} \left(g_0(t) + h \sum\nolimits_{j=1}^s a_{ij} U_j(t)\Omega_j \right).
	\]
Differentiating on both sides and letting $\dot{U}_i = U_i Y_i$, where $Y_i \in \sk (n)$ is a left-trivialized tangent vector, we obtain
	\[
		Y_i = d \mathbb{P}_{A_i}\left(g_0 \delta g_0 + h\sum\nolimits_{j=1}^s a_{ij} U_j Y_j\Omega_j\right),
	\]
	which can be rewritten as 
	\[
		Y_i - d \mathbb{P}_{A_i}\left(h \sum\nolimits_{j=1}^s a_{ij} U_j Y_j \Omega_j \right) = d\mathbb{P}_{A_i}(g_0 \delta g_0).
	\]
	Similar to the approach used for $X_{ik}$, we introduce a new map,
	\[
		\tilde{\varphi} \colon \sk (n) \longrightarrow \sk (n)^s, \quad \delta g_0 \longmapsto \left\{d \mathbb{P}_{A_i}(g_0 \delta g_0)\right\}_{i=1}^s,
	\]
	then $Y_i = (\tilde{\pi}_i \circ \tilde{\phi}^{-1}\circ \tilde{\varphi})(\delta g_0)$. The explicit expression for $\tilde{\varphi}^*$ can be written as (see section \ref{append_varphistar} for details),
	\begin{equation}
	\label{varphistar}
		\tilde{\varphi}^*(S_1, \ldots, S_s) = \asym \left(g_0^{T}\sum\nolimits_{i=1}^s d\mathbb{P}_{A_i}^*(S_i)\right).
	\end{equation}
	As such, $(\tilde{\pi}_i\circ \tilde{\phi}^{-1}\circ \tilde{\varphi})^*(S)$ can be computed as follows,
	\begin{subequations}
	\label{fff}
	\begin{align}[left = \empheqlbrace\,]
			S_j - \asym \left(h U_j^{T} \sum\nolimits_{l=1}^s a_{lj} d \mathbb{P}_{A_l}^*(S_l) \Omega_j^{T}\right) &= S \cdot \delta_{ij},	\qquad j = 1, 2, \ldots s,	\label{fff1}\\
			(\tilde{\pi}_i \circ \tilde{\phi}^{-1} \circ \tilde{\varphi})^*(S) &= \asym \left(g_0^{T} \sum\nolimits_{l=1}^s d \mathbb{P}_{A_l}^*(S_l) \right).	\label{fff2}
	\end{align}
	\end{subequations}
	Now, we compute $\frac{\partial \tilde{F}}{\partial g_0}$, which gives us \eqref{g0} explicitly (see section \ref{append_dFdg0} for details),
	\begin{equation}
	\label{ng0}
	 - p_0 = h \sum\limits_{i=1}^s b_i (\tilde{\pi}_i \circ \tilde{\phi}^{-1} \circ \tilde{\varphi})^* \left( \frac{\partial L}{\partial U}(U_i,\Omega_i) + \asym (U_i^{T} g_1 \Lambda \Omega_i^{T}) \right) + \asym (g_0^{T} g_1 \Lambda).
	\end{equation}
	
	 For equation \eqref{g1}, it is easy to show that
	\begin{equation}
	\label{ng1}
    	p_1 = \asym \left( g_1^{T} \left(g_0 + h\sum\nolimits_{i=1}^s b_i U_i \Omega_i \right) \Lambda^{T}\right).
	\end{equation}
	Combining \eqref{nomega}, \eqref{restrict}, \eqref{internal}, \eqref{ng0}, and \eqref{ng1}, we obtain a Lagrangian variational integrator on $SO(n)$,
	\begin{subequations}
	\label{vil}
	\begin{align}[left = \empheqlbrace\,]
		\begin{split}
			0 &= h \sum_{i=1}^s b_i (\tilde{\pi}_i \circ \tilde{\phi}^{-1}\circ \tilde{\psi}_k)^* \left( \frac{\partial L}{\partial U} (U_i, \Omega_i) + \asym (U_i^{T} g_1 \Lambda\Omega_i^{T}) \right)	\\
			& \qquad +  b_k \left( \frac{\partial L}{\partial \Omega}(U_k,\Omega_k) + \asym (U_k^{T} g_1 \Lambda) \right)
		\end{split}	\label{vil1}	\\
			g_1 &= \mathbb{P}\left(g_0 + h\sum\nolimits_{i=1}^s b_i U_i \Omega_i \right),	\label{vil2}		\\
			U_i &= \mathbb{P}\left(g_0 + h\sum\nolimits_{j=1}^s a_{ij} U_j \Omega_j \right),		\label{vil3}	\\
			 - p_0 &= h \sum\limits_{i=1}^s b_i (\tilde{\pi}_i \circ \tilde{\phi}^{-1} \circ \tilde{\varphi})^* \left( \frac{\partial L}{\partial U}(U_i,\Omega_i) + \asym (U_i^{T} g_1 \Lambda \Omega_i^{T}) \right) + \asym (g_0^{T} g_1 \Lambda),	\label{vil4}		\\
			 p_1 &= \asym \left( g_1^{T} \left(g_0 + h\sum\nolimits_{i=1}^s b_i U_i \Omega_i \right) \Lambda^{T}\right).	\label{vil5}	
	\end{align}
	\end{subequations}
	The integrator gives an update map $(g_0,p_0) \mapsto (g_1,p_1)$ on the cotangent bundle. In particular, one may solve for $(\{\Omega\}_{k=1}^s,g_1,\{U_i\}_{i=1}^s,\Lambda)$ simultaneously using equations \eqref{vil1}--\eqref{vil4}. Unfortunately, while \eqref{vil2}--\eqref{vil4} can be written in fixed-point form for the variables $g_1$, $\{U_i\}_{i=1}^s$, and $\Lambda$, \eqref{vil2} is implicit for $\{\Omega\}_{k=1}^s$. However, we can arrive at a fixed-point form for \eqref{vil1} on the Hamiltonian side if $L$ is hyperregular.

\section{Hamiltonian Variational Integrators on the Rotation Group $SO(n)$}
\label{sect:HI}
	It is often desirable to transform a numerical method from the Lagrangian side to the Hamiltonian side, which is possible if $L$ is hyperregular. The same mechanical system can be represented by either a Lagrangian or a Hamiltonian, and they are related by the Legendre transform. In Euclidean space, this gives
 	\[
		\xymatrix{
			(TQ, L)	\ar@<0.5ex>	[r]^{\mathbb{F}L} & (T^*Q,H)\ar@<0.5ex> [l]^{\mathbb{F}H},
 		}
 	\]
	and we have the following relationships,
	\[
		\frac{\partial L}{\partial \dot{q}}(q,\dot{q}) = p, 
		\quad \frac{\partial H}{\partial p}(q,p) = \dot{q},
		\quad \frac{\partial L}{\partial q}(q,\dot{q}) = -\frac{\partial H}{\partial q}(q,p).
	\]
	Given a Lie group $G$, a left-trivialized Lagrangian $L \colon G\times \mathfrak{g}\rightarrow\mathbb{R}$, and its corresponding Hamiltonian $H \colon G\times \mathfrak{g}^*\rightarrow\mathbb{R}$, it is easy to verify that similar relationships hold for these trivializations,
	\begin{equation}
	\label{transform}
		\frac{\partial L}{\partial \epsilon}(g,\epsilon) = \mu,
		\quad \frac{\partial H}{\partial \mu}(g,\mu) = \epsilon,
		\quad \frac{\partial L}{\partial g}(g,\epsilon) = -\frac{\partial H}{\partial g}(g,\mu).
	\end{equation}
	Using \eqref{transform} and denoting the corresponding internal cotangent vectors by $\{\mu_k\}_{k=1}^s$, \eqref{vil} can be transformed to the Hamiltonian form as
	\begin{subequations}
	\label{vih}
	\begin{align}[left = \empheqlbrace\,]
			\mu_k &= -\asym (U_k^{T} g_1 \Lambda) + h \sum\limits_{i=1}^s \frac{b_i}{b_k} (\tilde{\pi}_i \circ \tilde{\phi}^{-1} \circ \tilde{\psi}_k)^* \left( \frac{\partial H}{\partial U} (U_i, \mu_i) - \asym (U_i^{T} g_1 \Lambda \Omega_i^{T})\right),	\label{vih1}	\\
			g_1 &= \mathbb{P} \left(g_0 + h\sum\nolimits_{i=1}^sb_iU_i\Omega_i\right),	\label{vih2}	\\
			U_i &= \mathbb{P} \left(g_0 + h\sum\nolimits_{j=1}^s a_{ij} U_j \Omega_j\right),	\label{vih3}	\\
			\asym (g_0^{T} g_1 \Lambda) &= -p_0 + h \sum\limits_{i=1}^s b_i(\tilde{\pi}_i\circ \tilde{\phi}^{-1}\circ \tilde{\varphi})^* \left( \frac{\partial H}{\partial U}(U_i, \mu_i) - \asym (U_i^{T} g_1 \Lambda \Omega_i^{T}) \right),	\label{vih4}	\\
			 p_1 &= \asym \left( g_1^{T} \left(g_0 + h\sum\nolimits_{i=1}^s b_i U_i \Omega_i \right) \Lambda^{T}\right),	\label{vih5}	\\
			\Omega_i &= \frac{\partial H}{\partial \mu}(U_i,\mu_i).	\label{vih6}	
	\end{align}
	\end{subequations}
	In the above algorithm, $\Omega_i$ is given explicitly by \eqref{vih6} and only serves to reduce the redundancy in the computations because they are used numerous times in the other expressions. Similarly, $(\tilde{\pi}_i\circ \tilde{\phi}^{-1})^*$ shows up in both \eqref{vih1} and \eqref{vih4}, so one can save computational effort by reusing the shared solution to \eqref{f1} and \eqref{fff1}. Then, the first four equations can be solved simultaneously by fixed-point iterations, meaning the variables $(\{\mu_k\}_{k=1}^s,g_1,\{U_i\}_{i=1}^s,\Lambda)$ are updated concurrently in each iteration. Observe that the fixed-point form for $\Lambda$ in \eqref{vih4} requires solving a Lyapunov equation. Finally, $p_1$ is solved explicitly in \eqref{vih5} after solving for $(\{\mu_k\}_{k=1}^s,g_1,\{U_i\}_{i=1}^s,\Lambda)$. We shall call the integrators defined by \eqref{vih} the \textit{variational polar decomposition} method or VPD for short. 
	\subsection{Lie--Poisson Integrator by Reduction}
	\label{sect:HILiePoisson}
	We consider a $G$-invariant Hamiltonian system given by $H$ on the cotangent bundle $T^*G$. In this case, Hamilton's equations can be reduced to a Lie--Poisson system on $\mathfrak{g}^*$. If the Hamiltonian is hyperregular, then both the Lagrangian and the corresponding constrained Galerkin discrete Lagrangian $L_d(g_0,g_1)$ will be $G$-invariant. As such, \eqref{vi} naturally reduces to yield a Lie--Poisson integrator (see Appendix \ref{append_LPReductions}). We only consider the reduction on the Hamiltonian side due to the nature of the constrained Galerkin methods, which give an update map on the cotangent bundle. 
	
The discrete Lagrangian we have constructed becomes
	\[
		L_d(g_0,g_1) = \ext_{\{\Omega_i\}_{i=1}^s} \sum\limits_{i=1}^s b_i \bm{l}(\Omega_i),
	\]
	where 
	\begin{align*}
		U_i &= \mathbb{P}\left(g_0 + h \sum\nolimits_{j=1}^s a_{ij} U_j \Omega_j\right),	\\
		 g_1 &= \mathbb{P} \left(g_0 + h \sum\nolimits_{i=1}^s b_i U_i \Omega_i\right),
	\end{align*}
and $\bm{l} \colon \mathfrak{g} \to \mathbb{R}$ is the reduced Lagrangian. It is easy to verify that our system is $G$-invariant, meaning
	\[
		L_d(g\cdot g_0,g\cdot g_1) = L_d(g_0,g_1) = \ext_{\{\Omega_i\}_{i=1}^s} \sum\limits_{i=1}^s b_i \bm{l}(\Omega_i),
		\]
where
	\begin{align*}
		g \cdot U_i &= \mathbb{P} \left((g \cdot g_0) + h \sum\nolimits_{j=1}^s a_{ij} (g\cdot U_j \Omega_j)\right),\\
		g \cdot g_1 &= \mathbb{P} \left((g \cdot g_0) + h \sum\nolimits_{i=1}^s b_i (g \cdot U_i \Omega_i) \right).
\end{align*}
	Therefore, the variational integrator \eqref{vih} can theoretically be reduced to a Lie--Poisson integrator. By letting $g_0^{T}g_1=f_0$ and $U_i^{T}g_1=\Theta_i$, \eqref{vih} simplifies to
	\begin{subequations}
	\label{lp}
	\begin{align}[left = \empheqlbrace\,]
			\mu_k &= -\asym (\Theta_k \Lambda) - h \sum\limits_{i=1}^s \frac{b_i}{b_k} (\tilde{\pi}_i \circ \tilde{\phi}^{-1} \circ \tilde{\psi}_k)^*(\asym (\Theta_i \Lambda \Omega_i^{T})),	\label{lp1}	\\
			0 &= \asym \left( f_0^{T} + h \sum\nolimits_{i=1}^s b_i \Theta_i^{T}\Omega_i \right), \label{lp2}\\
			\Theta_i^{T} &= \mathbb{P} \left( f_0^{T} + h \sum\nolimits_{j=1}^s a_{ij} \Theta_j^{T} \Omega_j \right), \label{lp3}\\
			\asym (f_0 \Lambda) &= -p_0 - h \sum\limits_{i=1}^s b_i (\tilde{\pi}_i \circ \tilde{\phi}^{-1} \circ \tilde{\varphi})^*(\asym (\Theta_i \Lambda \Omega_i^{T})),	\label{lp4}\\
			p_1 &= \asym \left( \left(f_0^{T} + h \sum\nolimits_{i=1}^s b_i \Theta_i^{T} \Omega_i \right)\Lambda^{T} \right),	\label{lp5}\\
			\Omega_i &= \frac{\partial \bm{h}}{\partial \mu}(\mu_i), \label{lp6}
	\end{align}
	\end{subequations}
	where $\bm{h}$ is the reduced Hamiltonian. Multiplying by $g_1^{T}$ on both sides of $g_1 = \mathbb{P} \left(g_0 + h\sum\nolimits_{i=1}^s b_i U_i \Omega_i \right)$ yields 
	\[
		I = \mathbb{P} \left(f_0^{T} + h \sum\nolimits_{i=1}^s b_i \Theta_i^{T} \Omega_i \right).
	\]
	Suppose that $h$ is small and $g_0$ and $g_1$ are close, then $\left(f_0^{T} + h\sum\nolimits_{i=1}^s b_i \Theta_i^{T} \Omega_i \right)$ is in the neighborhood of $I$. By Lemma~\ref{null}, this is equivalent to
	\[
		\asym \left(f_0^{T} + h \sum\nolimits_{i=1}^s b_i \Theta_i^{T} \Omega_i \right) =0,
	\]
	which can be regarded as the fixed-point equation for $f_0$. The first four equations can be solved using fixed-point iteration for the variables $(\{\mu_k\}_{k=1}^s,f_0,\{\Theta_i\}_{i=1}^s,\Lambda)$ as in our previous discussions. Then, $p_1$ can be calculated explicitly. 
	
	In the above algorithm, we also need to figure out the reduced version of $(\tilde{\pi}_i \circ \tilde{\phi}^{-1} \circ \tilde{\psi}_k)^*$ and $(\tilde{\pi}_i \circ \tilde{\phi}^{-1} \circ \tilde{\varphi})^*$. Note that \eqref{f1} and \eqref{f2} involve $U_j^{T}d \mathbb{P}_{A_i}^*$; in addition, $g_0, g_1, U_i$, and $A_i = g_0 + h \sum\nolimits_{j=1}^s a_{ij} U_j \Omega_j = U_i P_i$ are reduced, so we need a reduced version of $U_j^{T} d \mathbb{P}_{A_i}^*$ as well. Multiplying $A_i$ on the left by $g_1^T$, we obtain
	\[
		g_1^{T} A_i = f_0^{T} + h \sum\limits_{j=1}^s a_{ij} \Theta_j^{T} \Omega_j = (g_1^{T} U_i) P_i.
	\]
Then, $(g_1^{T} U_i) P_i$ is the polar decomposition of $f_0^{T} + h \sum\nolimits_{j=1}^s a_{ij} \Theta_j^{T} \Omega_j$, and for $S \in \sk (n)$,
	\begin{align*}
		U_j^{T} d \mathbb{P}_{A_i}^*(S) &= U_j^{T} \cdot U_i \lyap (P_i,S^{T}) = \Theta_j g_1^{T} U_i \lyap (P_i,S^{T}),	\\
		&= \Theta_j \mathbb{P} \left(f_0^{T} + h \sum\nolimits_{j=1}^s a_{ij} \Theta_j^{T} \Omega_j \right) \lyap (P_i,S^{T}).
	\end{align*}
	This is the reduced version of $U_j^{T} d \mathbb{P}_{A_i}^*(S)$, and so $(\tilde{\pi}_i \circ \tilde{\phi}^{-1} \circ \tilde{\psi}_k)^*(S)$ can be computed as follows,
	\begin{subequations}
	\label{nf}
	\begin{align}[left = \empheqlbrace\,]
		S_j - \asym \left(h \Theta_j \sum\nolimits_{l=1}^s a_{lj} d \mathbb{P}_{A_l}^* (S_l) \Omega_j^{T} \right) = S \cdot \delta_{ij},	\qquad j = 1, 2 \ldots s,\label{nf1}	\\
		(\tilde{\pi}_i \circ \tilde{\phi}^{-1} \circ \tilde{\psi}_k)^*(S) = \asym \left( \Theta_k \sum\nolimits_{l=1}^s a_{lk} d \mathbb{P}_{A_l}^*(S_l)\right), 	\label{nf2}
	\end{align}	
	\end{subequations}
	where $\{A_i\}$ are redefined to be $A_i = f_0^{T} + h \sum\nolimits_{j=1}^s a_{ij} \Theta_j^{T} \Omega_j$. For $(\tilde{\pi}_i \circ \tilde{\phi}^{-1} \circ \tilde{\varphi})^*$, we need to compute $g_0^{T} d \mathbb{P}_{A_i}^*(S)$ for some $S \in \sk (n)$, which is given by
	\begin{align*}
		g_0^{T} d \mathbb{P}_{A_i}^*(S) 
		&= g_0^{T} \cdot U_i \lyap (P_i,S^{T}) = f_0 g_1^{T} U_i \lyap (P_i, S^{T}) 	\\
		&= f_0 \mathbb{P} \left( f_0^{T} + h \sum\nolimits_{j=1}^s a_{ij} \Theta_j^{T} \Omega_j\right) \lyap (P_i,S^{T}).
	\end{align*}
	Hence, we have
	\begin{subequations}
	\label{nnf}
	\begin{align}[left = \empheqlbrace\,]
		S_j - \asym \left(h \Theta_j \sum\nolimits_{l=1}^s a_{lj} d \mathbb{P}_{A_l}^* (S_l) \Omega_j^{T} \right) = S \cdot \delta_{ij},	\qquad j = 1, 2 \ldots s,\label{nnf1}	\\
		(\tilde{\pi}_i \circ \tilde{\phi}^{-1} \circ \tilde{\varphi})^*(S) = \asym \left( f_0 \sum\nolimits_{l=1}^s a_{lk} d \mathbb{P}_{A_l}^*(S_l)\right). 	\label{nnf2}
	\end{align}	
	\end{subequations}

\section{Numerical Experiments}
\label{sect:NumericalExp}
	
	We test our methods and compare them to the variational Runge--Kutta--Munthe-Kaas (VRKMK) methods from \citet{Geir2014} on the dipole on a stick problem that they considered. In particular, our configuration space is $SO(3)$, and its Lie algebra is identified with $\mathbb{R}^3$. We shall only recall the mathematical expressions here, so for a thorough description of the system one should refer to the reference above. The right-trivialized and left-trivialized Hamiltonians, $H^R,H^L \colon SO(3) \times \mathbb{R}^3 \to \mathbb{R}$, can be written as
	\begin{align}
		H^R(g,p) 
			&= \frac{1}{2} p^T g J^{-1} g^T p + U(g),	\label{RightH}	\\
		H^L(g,\tilde{p}) 
			&= \frac{1}{2} \tilde{p}^T J^{-1} \tilde{p} + U(g),	\label{LeftH}
	\end{align}
	where $U(g) =  me_3^T g e_3 + q\beta\left( \lVert gy_{+}^{0} - z \rVert_2^{-1} - \lVert gy_{-}^{0} - z \rVert_2^{-1} \right)$. Note that $J = m \, \text{diag}(1+\alpha^2,1,\alpha^2)$, with $m = 1$ and $\alpha = 0.1$. The constant vectors are $y_{\pm}^0 = (0,\pm \alpha, -1)^T$ and  $z=(0,0,-3/2)^T$. Lastly, $q = \beta = 1$, and $\|\cdot\|_2$ is our usual Euclidean norm.
	
	Both forms of the Hamiltonian are written here because while our method was developed for left-trivialized systems using Hamilton's principle, their method was developed for right-trivialized systems using the Hamilton--Pontryagin principle. As a result, both discretizations yield symplectic variational integrators for the Hamiltonian with their corresponding choice of trivialization. In particular, we have
	\begin{equation}
	\label{eqn:LRTrivializedIG}
		\tilde{p} = g^Tp
	\end{equation}
	as a relationship between the dual elements of the corresponding trivialization: $p$ is the dual representation of $\xi \in \mathfrak{g}$ for the right-trivialization $\dot{g} = \xi g$, and $\tilde{p}$ is the dual representation of $\tilde{\xi} \in \mathfrak{g}$ for the left-trivialization $\dot{g} = g\tilde{\xi}$. Consequently, we note that the left-trivialized cotangent vector $\frac{\partial H}{\partial g}$ in VPD \eqref{vih} is computed as
	\[
		\frac{\partial H}{\partial g}(g,\tilde{p}) = 
			\asym\left( g^{-1} \nabla_g H^L(g,\tilde{p}) \right),
	\]
	where $\nabla_g H^L(g,\tilde{p})$ is the matrix derivative of $H^L(g,\tilde{p})$ with respect to $g$. On the other hand, the right-trivialized cotangent vector is computed as
	\[
		\frac{\partial H}{\partial g}(g,p) = 
			\asym\left(\nabla_g H^R(g,p)g^{-1} \right),
	\]
	when implementing the VRKMK method. 
	
	For our tests, we also have the same initial data from \citep{Geir2014},
	\begin{align*}
		g(0) 
			&=
			\begin{pmatrix}
				1	&	0	&	0	\\
				0	&	0	&	-1	\\
				0	&	1	&	0
			\end{pmatrix},	\\
		p(0)
			&= g(0) J g(0)^T e_2,
	\end{align*}
	 for the VRKMK methods, and so \eqref{eqn:LRTrivializedIG} gives $\tilde{p}(0)$ to complete the initial data for the VPD methods. Since both methods involve fixed-point iteration,  we terminate the processes when the norm $\|\cdot\|_2$ between the current and previous iteration is less than $10^{-15}$ for each variable. In particular, the norm for vectors is the Euclidean norm and for matrices it is the induced matrix norm. Lastly, we ran these implementations in Wolfram Mathematica 12 on a personal computer with the following specifications: Operating system: Windows 10; CPU: 3.7GHz AMD Ryzen 5 5600X; RAM: G.Skill Trident Z Neo Series 32Gb DDR4 3600 C16; Storage: 1TB Rocket Nvme PCIe 4.0 writing/reading up to 5000/4400 MB/s.
	\subsection{Order Tests}
	\label{sect:NE_OrderTests}
\begin{table*}[t]
	\begin{subtable}[b]{.45\linewidth}
		\centering
	\begin{tabular}{c|c}
		$\dfrac{1}{2}$ & $\dfrac{1}{2}$ \\[2ex] \hline \\[-2ex]
		               & 1             
	\end{tabular}
	\caption{2\textsuperscript{nd} order Gauss--Legendre method}
	\label{tab:2ndGL}
	\end{subtable}
	\vspace{5mm}
	\begin{subtable}[b]{.45\linewidth}
		\centering
	\begin{tabular}{c|ccc}
		$0$            & $0$            & $0$            & $0$            \\[1ex]
		$\dfrac{1}{2}$ & $\dfrac{1}{2}$ & $0$            & $0$            \\[2ex]
		$1$            & $-1$           & $2$            & $0$            \\[1ex] \hline \\[-2ex]
		               & $\dfrac{1}{6}$ & $\dfrac{2}{3}$ & $\dfrac{1}{6}$
	\end{tabular}
	\caption{3\textsuperscript{rd} order Runge--Kutta method}
	\label{tab:3rdRK}
	\end{subtable}
	\begin{subtable}[b]{.45\linewidth}
		\centering
	\begin{tabular}{c|cc}
		$\dfrac{1}{2}-\dfrac{\sqrt{3}}{6}$ & $\dfrac{1}{4}$                     & $\dfrac{1}{4}-\dfrac{\sqrt{3}}{6}$ \\[2ex]
		$\dfrac{1}{2}+\dfrac{\sqrt{3}}{6}$ & $\dfrac{1}{4}+\dfrac{\sqrt{3}}{6}$ & $\dfrac{1}{4}$                     \\[2ex] \hline \\[-2ex]
		                                   & $\dfrac{1}{2}$                  	& $\dfrac{1}{2}$                    
\end{tabular}
	\caption{4\textsuperscript{th} order Gauss--Legendre method}
	\label{tab:4thGL}
	\end{subtable}
	\begin{subtable}[b]{.45\linewidth}
		\centering
	\begin{tabular}{c|ccc}
		$\dfrac{1}{2}-\dfrac{\sqrt{15}}{10}$ & $\dfrac{5}{36}$                       & $\dfrac{2}{9}-\dfrac{\sqrt{15}}{15}$ & $\dfrac{5}{36}-\dfrac{\sqrt{15}}{30}$ \\[2ex]
		$\dfrac{1}{2}$                       & $\dfrac{5}{36}+\dfrac{\sqrt{15}}{24}$ & $\dfrac{2}{9}$                       & $\dfrac{5}{36}-\dfrac{\sqrt{15}}{24}$ \\[2ex]
		$\dfrac{1}{2}+\dfrac{\sqrt{15}}{10}$ & $\dfrac{5}{36}+\dfrac{\sqrt{15}}{30}$ & $\dfrac{2}{9}+\dfrac{\sqrt{15}}{15}$ & $\dfrac{5}{36}$                       \\[2ex] \hline \\[-2ex]
		                                     & $\dfrac{5}{18}$                       & $\dfrac{4}{9}$                       & $\dfrac{5}{18}$                      
	\end{tabular}
	\caption{6\textsuperscript{th} order Gauss--Legendre method}
	\label{tab:6thGL}
	\end{subtable}
	\caption{Butcher tableaux for the comparison tests}
	\label{table:Butchertableaux}
\end{table*}

\begin{figure}[t]
	\begin{subfigure}[b]{.45\textwidth}
		\includegraphics{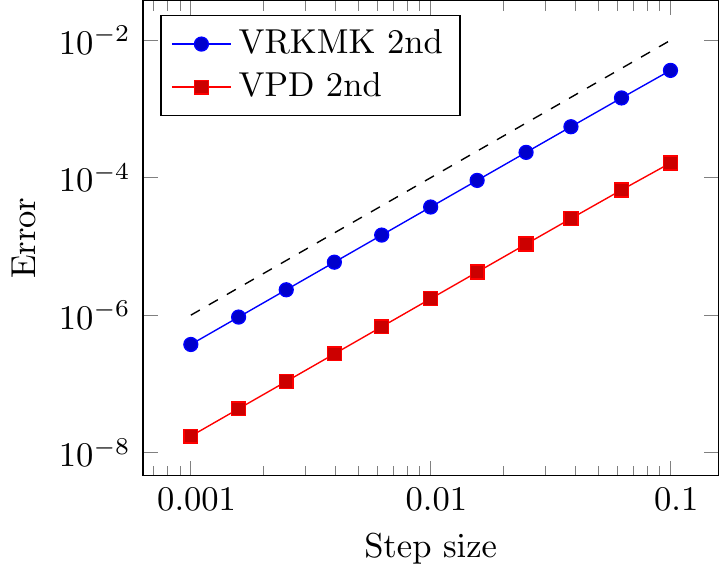}
	\end{subfigure}
	\begin{subfigure}[b]{.45\textwidth}
		\includegraphics{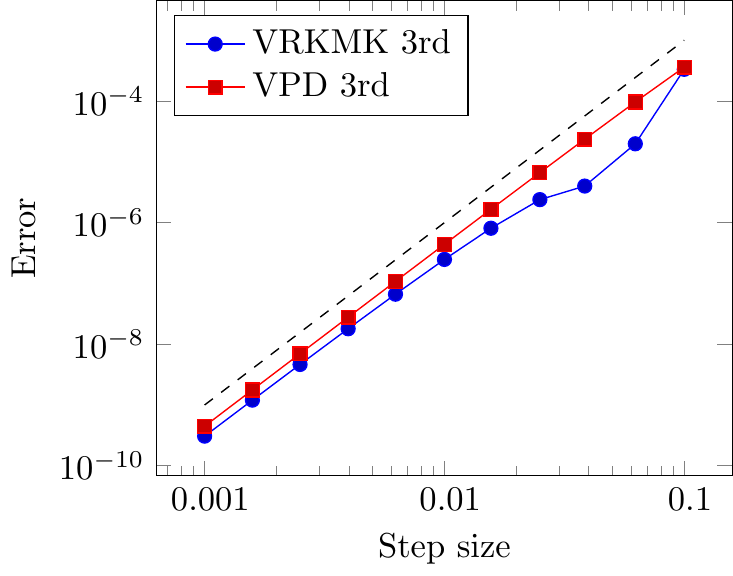}
	\end{subfigure}
	\begin{subfigure}[b]{.45\textwidth}
		\includegraphics{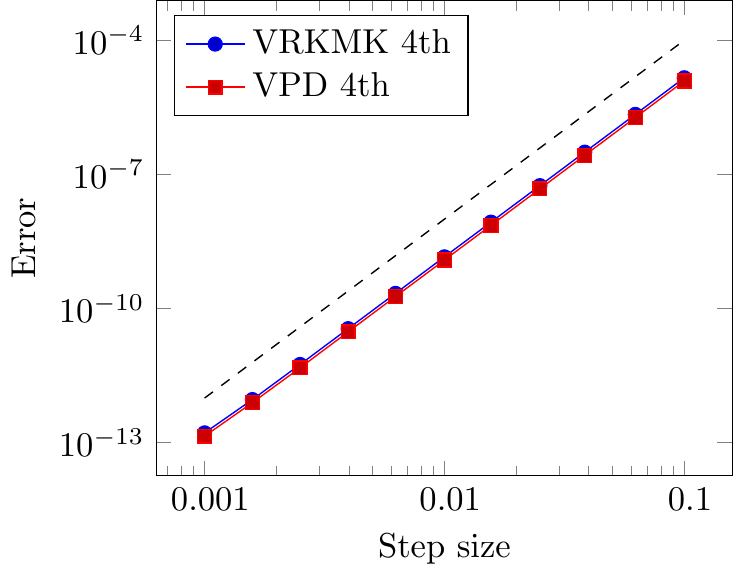}
	\end{subfigure}
	\begin{subfigure}[b]{.45\textwidth}
		\includegraphics{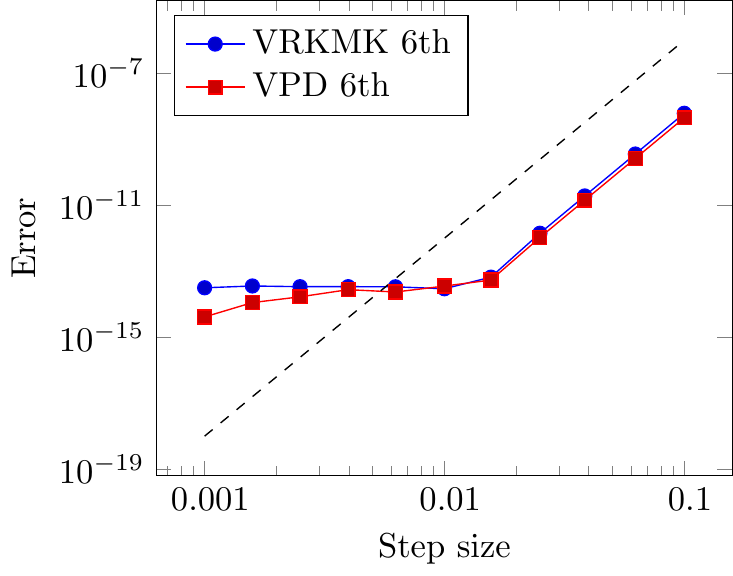}
	\end{subfigure}
	\caption{Order comparison plots between VRKMK and VPD methods: The black-dashed lines are references for the corresponding orders.}
	\label{fig:orderComp}
\end{figure}

	\begin{figure}[p]
	\begin{subfigure}[b]{1.\textwidth}
		\includegraphics{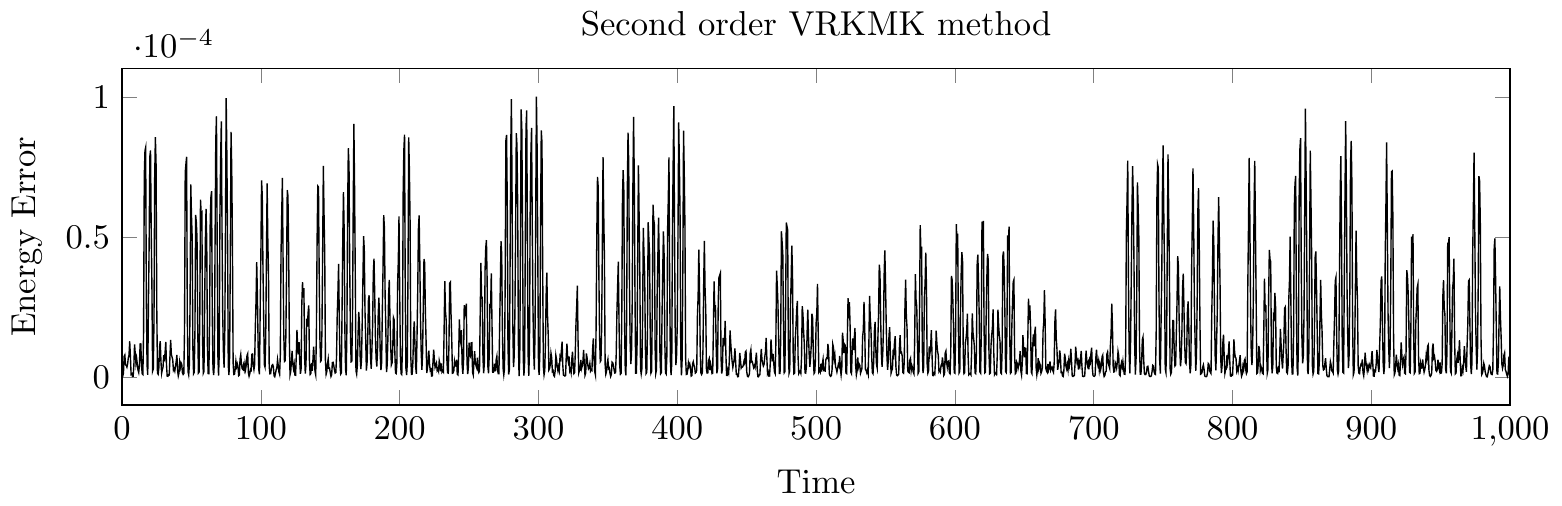}
	\end{subfigure}
	\begin{subfigure}[b]{1.\textwidth}
		\includegraphics{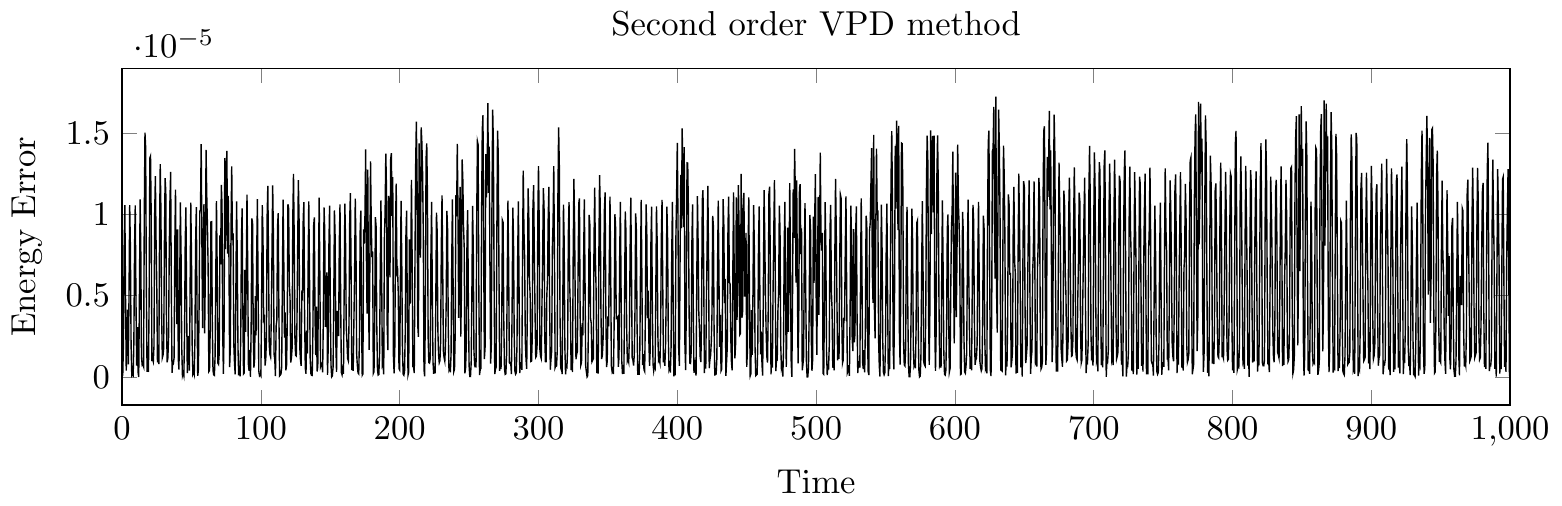}
	\end{subfigure}
	\caption{Long-term energy error for second order VRKMK/VPD methods.}
	\label{fig:EnergyErr2}
\end{figure}
\begin{figure}[p]
	\begin{subfigure}[b]{1.\textwidth}
		\includegraphics{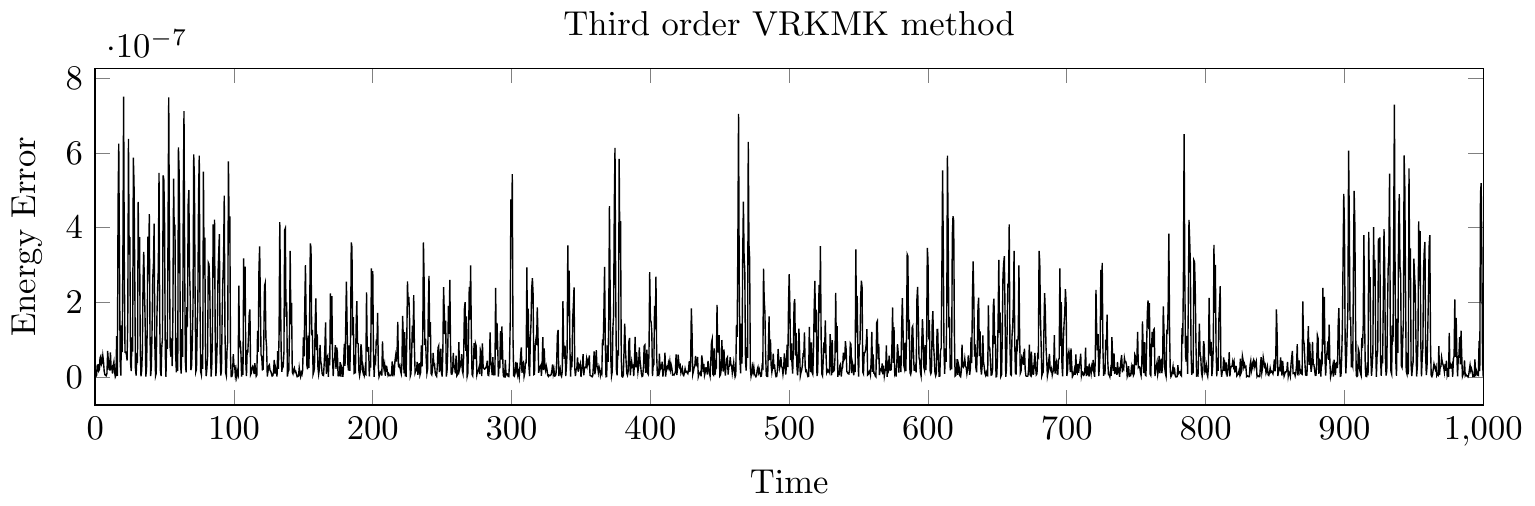}
	\end{subfigure}
	\begin{subfigure}[b]{1.\textwidth}
		\includegraphics{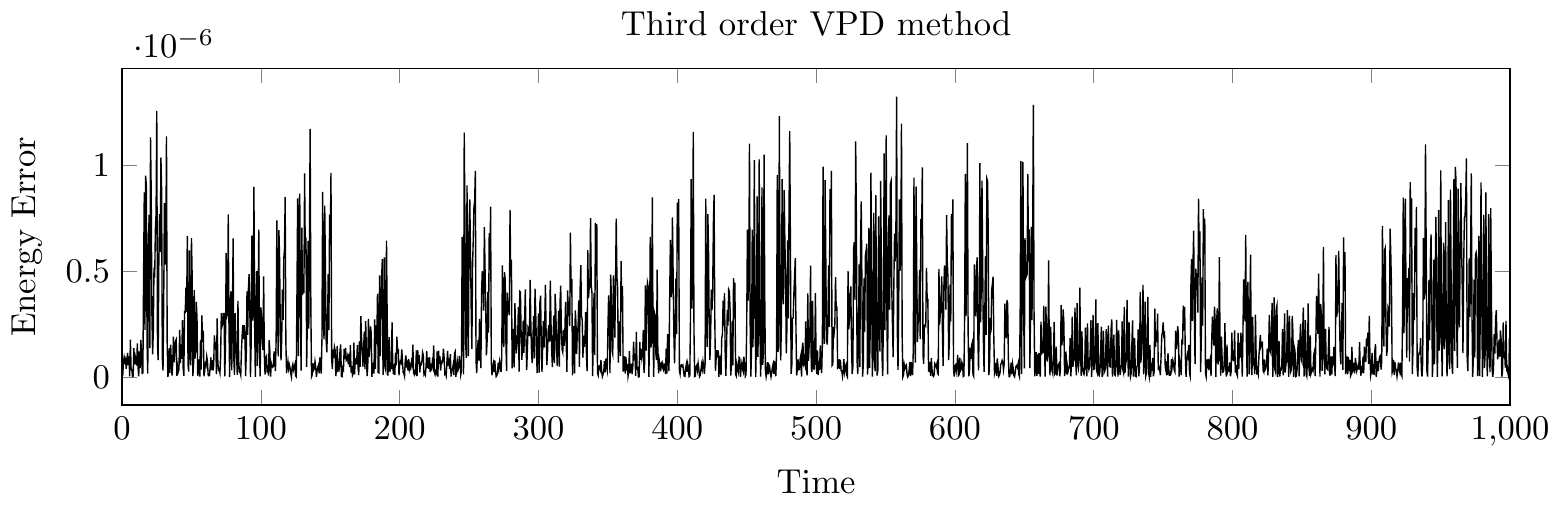}
	\end{subfigure}
	\caption{Long-term energy error for third order VRKMK/VPD methods.}
	\label{fig:EnergyErr3}
\end{figure}
\begin{figure}[p]
	\begin{subfigure}[b]{1.\textwidth}
		\includegraphics{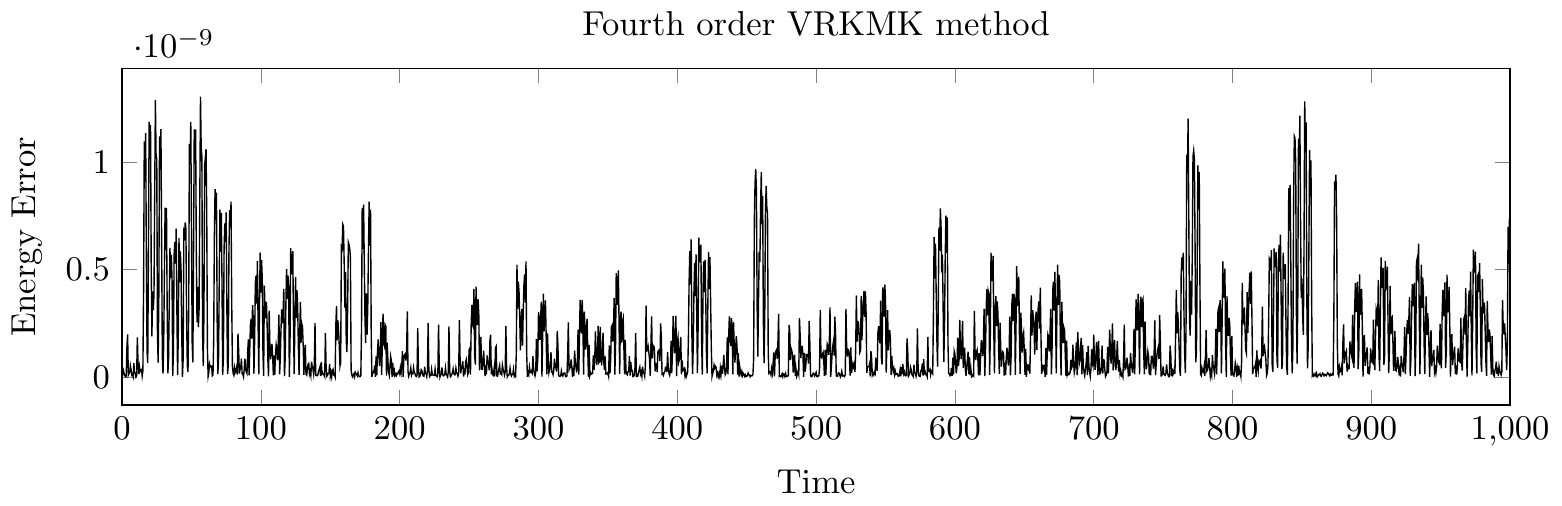}
	\end{subfigure}
	\begin{subfigure}[b]{1.\textwidth}
		\includegraphics{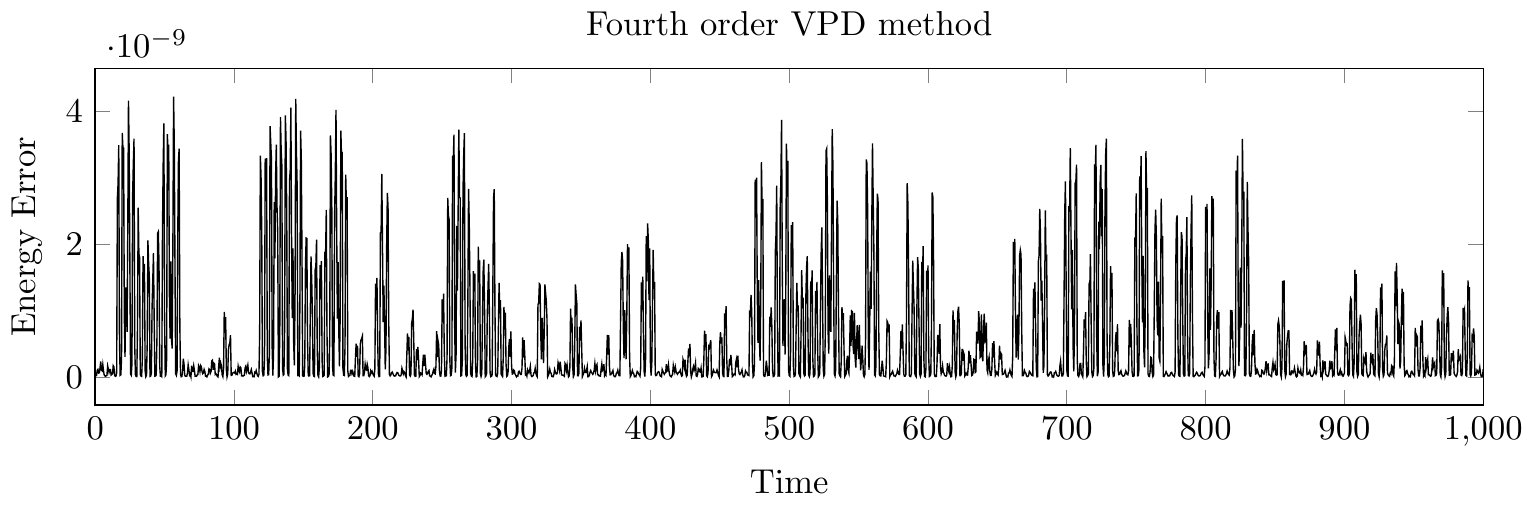}
	\end{subfigure}
	\caption{Long-term energy error for fourth order VRKMK/VPD methods.}
	\label{fig:EnergyErr4}
\end{figure}
\begin{figure}[p]
	\begin{subfigure}[b]{1.\textwidth}
		\includegraphics{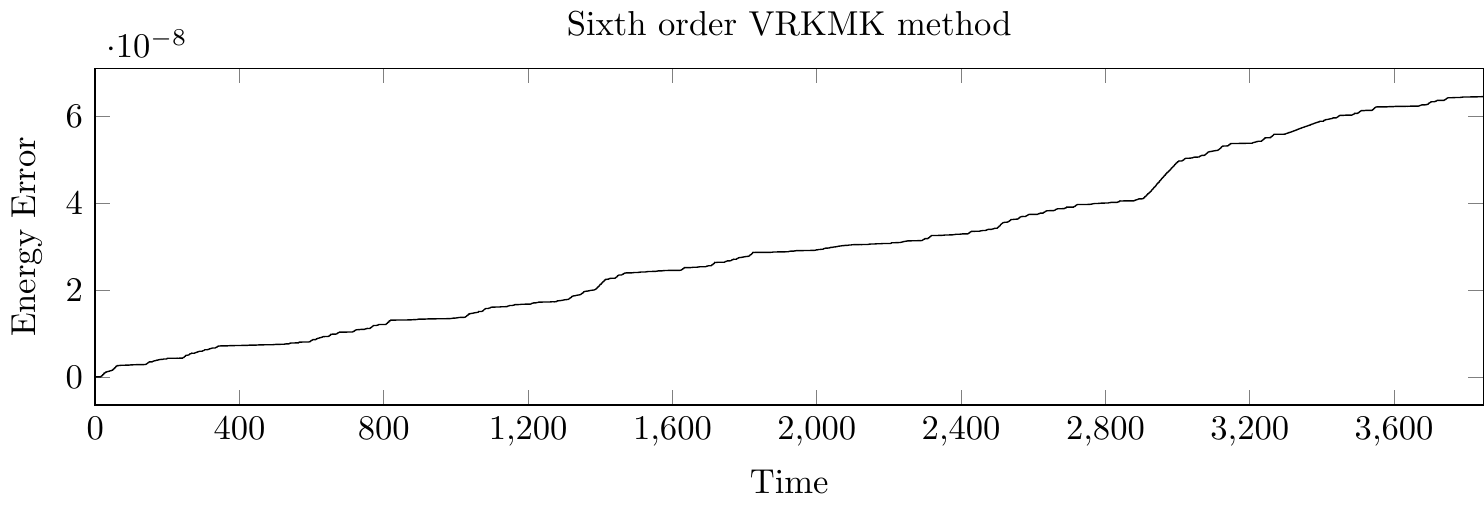}
	\end{subfigure}
	\begin{subfigure}[b]{1.\textwidth}
		\includegraphics{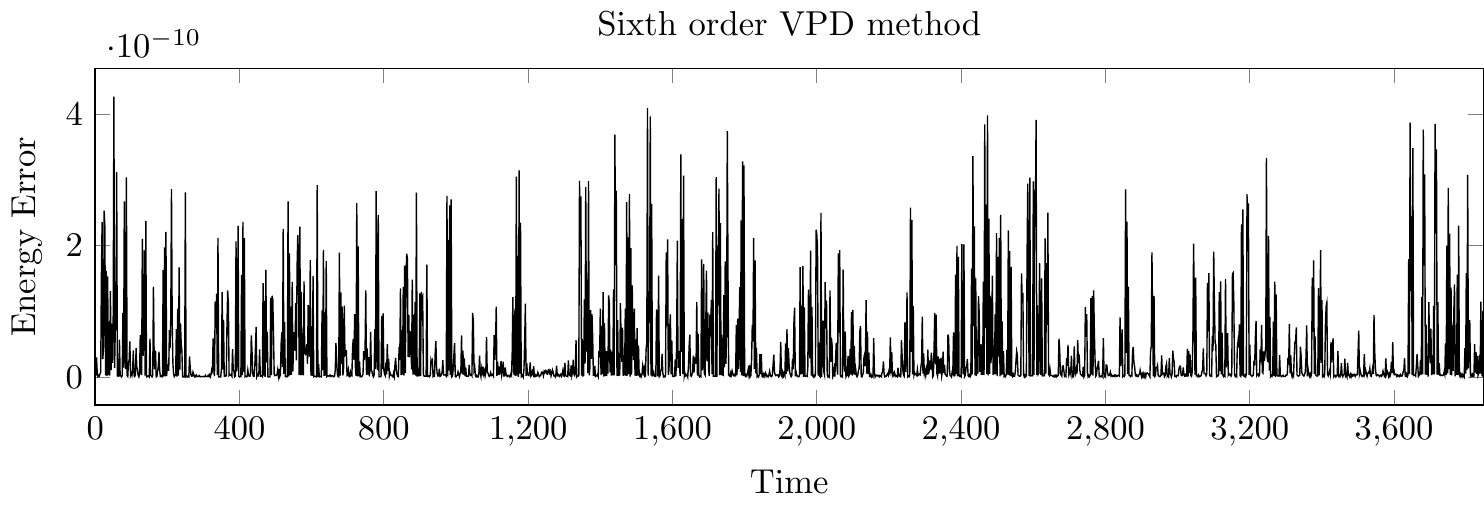}
	\end{subfigure}
	\caption{Long-term energy error for sixth order VRKMK/VPD methods.}
	\label{fig:EnergyErr6}
\end{figure}

\begin{figure}[h]
	\begin{subfigure}[b]{1.\textwidth}
		\includegraphics{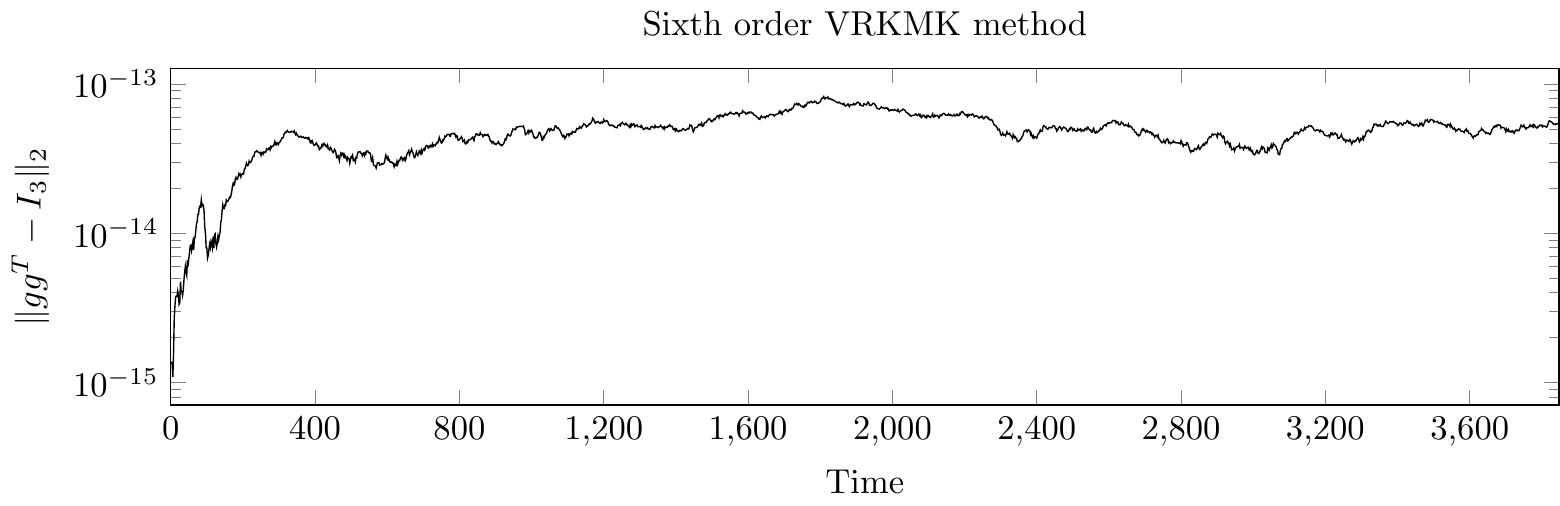}
	\end{subfigure}
	\begin{subfigure}[b]{1.\textwidth}
		\includegraphics{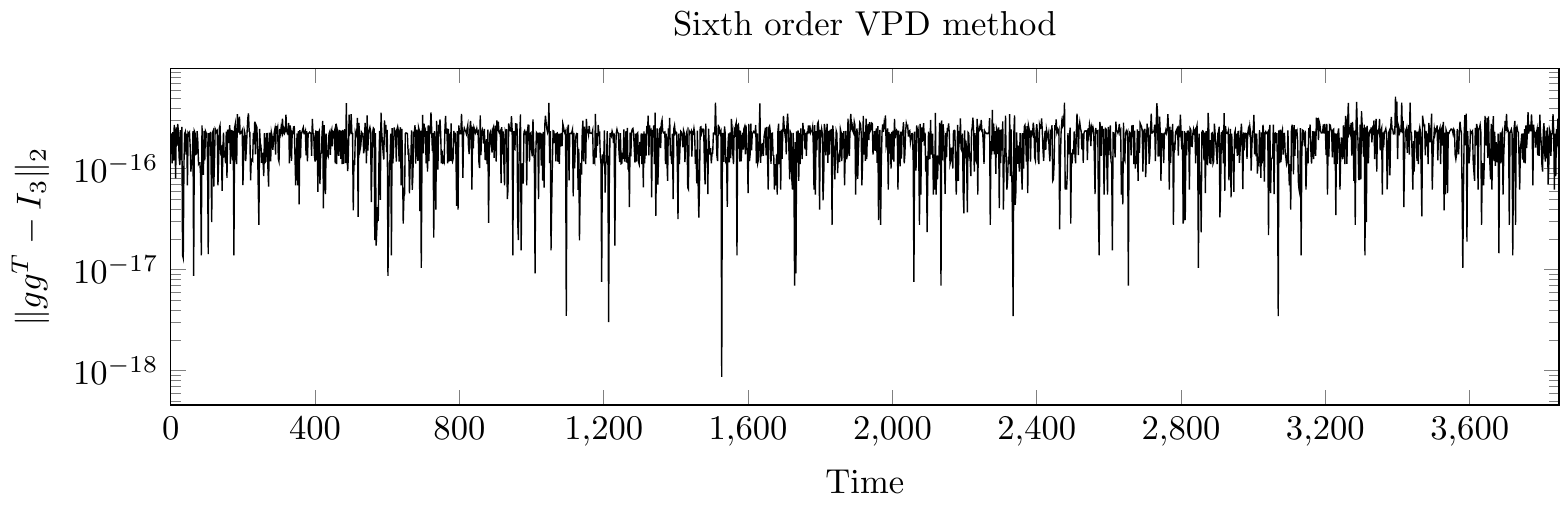}
	\end{subfigure}
	\caption{Orthogonality error for sixth order VRKMK/VPD methods.}
	\label{fig:OrthoErr6}
\end{figure}
\begin{figure}[h]
	\begin{subfigure}[b]{.475\textwidth}
		\includegraphics{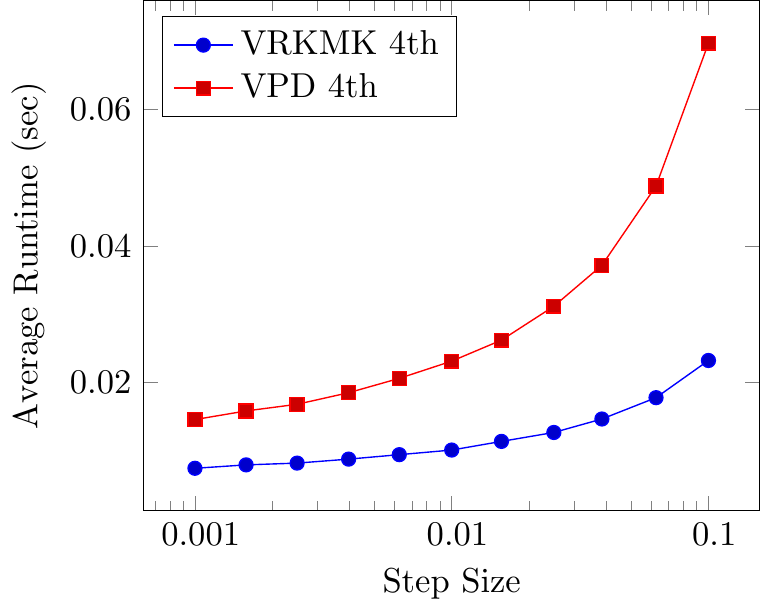}
	\end{subfigure}
	\begin{subfigure}[b]{.475\textwidth}
		\includegraphics{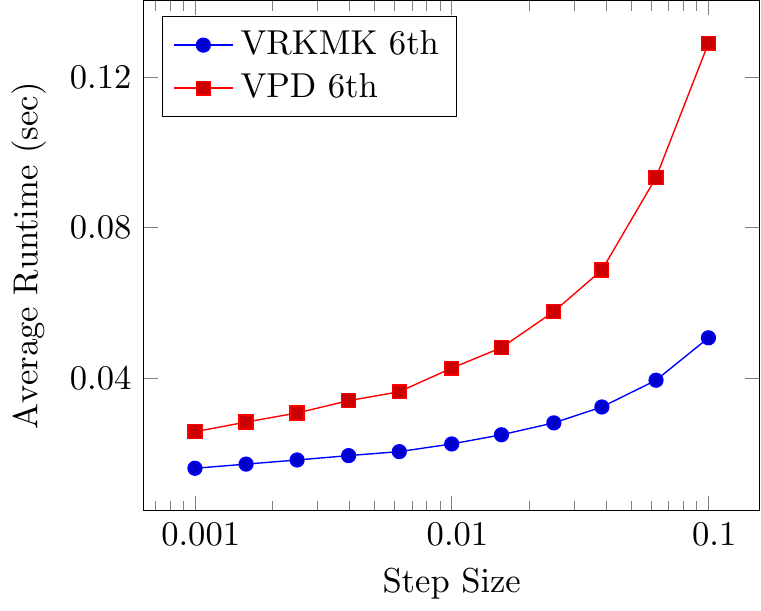}
	\end{subfigure}
	\caption{Run-time comparison for fourth and sixth order VRKMK/VPD methods.}
	\label{fig:Runtime46}
\end{figure}

	We run both methods based on the one-, two-, and three-stage Gauss--Legendre methods and a third-order Runge--Kutta method for comparisons, and these methods are shown as Butcher tableaux in Table \ref{table:Butchertableaux}. 
	
	We compute the errors in $(g(0.5),p(0.5))$ from VRKMK and errors in $(g(0.5),\tilde{p}(0.5))$ from VPD with respect to a reference solution and these errors are shown in Figure \ref{fig:orderComp}. The reference solution was calculated using the sixth-order VPD method with step size $h = 0.001$. The errors for VRKMK is computed as $\|g(0.5)^T p(0.5) - \tilde{p}_{\text{ref}}\|_2 + \| g(0.5) - g_{\text{ref}} \|_2 $ while the errors for VPD is computed as $\|\tilde{p}(0.5) - \tilde{p}_{\text{ref}}\|_2 + \| g(0.5) - g_{\text{ref}} \|_2 $. The black-dashed lines are reference lines for the appropriate orders: In the sixth order comparison plot,  the errors from the fixed-point iteration are dominating for smaller step sizes, so the theoretical order is not observed in those regimes. Otherwise, both methods achieve their theoretical orders, and they are quite comparable, noting that VPD exhibits a noticeable, smaller error constant in the second order method.
	
	\subsection{Long-Term Behaviors}
	\label{sect:NE_LTBehaviors}
	The long-term energy behaviors for both methods are presented in Figures \ref{fig:EnergyErr2} -- \ref{fig:EnergyErr6}. For second, third, and fourth order, we fixed the step size $h = 0.01$ and ran $10^5$ integration steps to observe the energy errors. For second order, the energy errors have magnitude orders of $10^{-4}$ for VRKMK and $10^{-5}$ for VPD; for third order, they have magnitude orders of $10^{-7}$ for VRKMK and $10^{-6}$ for VPD; for fourth order, they have magnitude orders of $10^{-9}$ for both methods. 
	
	For the sixth order, we consider step size $h = 1/26$ to avoid the regime where the errors in the fixed-point iteration are dominating. This step size corresponds to the third point from the right in the order comparison plots. We also run $10^5$ integration steps and observe the energy errors in Figure \ref{fig:EnergyErr6}. The energy error for the VPD method is stable with an order of magnitude of $10^{-10}$. On the other hand, the energy error in VRKMK exhibits a slow increase over the integrating time span. We investigated this phenomenon and attribute this to the framework of each method. In VPD, $g_1$ in \eqref{vih2} and the internal points $U_i$ in \eqref{vih3} are updated in each integration step via polar decomposition. As a result, the Lie group structure is always preserved up to machine precision for both the configuration space elements and internal points. However, in VRKMK, both $g_1$ and the internal points $U_i$ are updated by the left action of $SO(n)$. This left action is a matrix multiplication which is performed to machine precision, but with a sufficient number of multiplications, the round-off error will still accumulate. Consequently, the Lie group structure is not as well-preserved in comparison, and this is illustrated in the comparison of orthogonality errors $\|gg^T - I_3\|_2$ in Figure \ref{fig:OrthoErr6}. 
	
	We also observed that projecting the sixth-order VRKMK method onto the rotation group using the polar decomposition at each timestep does not recover the near energy preservation typical of symplectic methods. Presumably, this is because the projection subtly compromises the symplecticity of the method. This is consistent with the observation made in \cite{LeLeMc2007b}, that both the Lie group structure and symplecticity needs to be preserved for the methods to exhibit near energy preservation.
	\subsection{Runtime Comparison}
	\label{sect:NE_Runtime}
	We also have data for runtime comparison of VRKMK and VPD for the two-stage and three-stage Gauss--Legendre methods in Figure \ref{fig:Runtime46}. For each step size $h$, we recorded the runtime of each method in seconds and repeated the execution 256 times to compute the average runtime. The runtime for VPD methods is on average longer than for VRKMK methods. This may be due in parts to the Lyapunov solutions and multiple layers of fixed-point iterations required in the polar decomposition in \eqref{newtoniterate} and $(\tilde{\pi}_i \circ \tilde{\phi}^{-1})^*$ in the same equations \eqref{f1} and \eqref{fff1}.
	
\section{Conclusion}
\label{sect:Conclusion}
	By applying the polar decomposition in a constrained Galerkin method, we constructed high-order Lie group symplectic integrators on $SO(n)$, which we refer to as the variational polar decomposition method. Furthermore, the integrator can be reduced to a Lie--Poisson integrator when the underlying mechanical system is group-invariant. We performed extensive numerical experiments comparing our proposed method to the variational Runge--Kutta--Munthe-Kaas method from \citep{Geir2014}. These comparisons provide insights into each method and highlight an advantage of our proposed method, which is the preservation of Lie group structure up to machine precision for both the configurations and internal points. This appears to be important for high-order methods to exhibit the near energy preservation that one expects for symplectic integrators when applied to dynamics on Lie groups.
	
	For future work, it would be interesting to explore the generalization of the proposed method to symmetric spaces, by applying the generalized polar decomposition~\cite{munthe2001generalized}. This may be of particular interest in the construction of accelerated optimization algorithms on symmetric spaces, following the use of time-adaptive variational integrators for symplectic optimization~\cite{duruisseaux2020adaptive} based on the discretization of Bregman Lagrangian and Hamiltonian flows~\cite{WiWiJo16,Duruisseaux2021Riemannian}. Examples of symmetric spaces include the space of Lorentzian metrics, the space of symmetric positive-definite matrices, and the Grassmannian.

\section*{Acknowledgements}
This research has been supported in part by NSF under grants DMS-1010687, CMMI-1029445, DMS-1065972, CMMI-1334759, DMS-1411792, DMS-1345013, DMS-1813635, CCF-2112665, by AFOSR under grant FA9550-18-1-0288, and by the DoD under grant HQ00342010023 (Newton Award for Transformative Ideas during the COVID-19 Pandemic). The authors would also like to thank Geir Bogfjellmo and H\r{a}kon Marthinsen for sharing the code for their VRKMR method from \cite{Geir2014}, which we used in the numerical comparisons with our method.

\appendix

\section{Constrained Galerkin Methods}
\label{append_constrainedGalerkin}

Our Galerkin variational integrator will involve a discrete Lagrangian that differs from the classical construction in \cite{MaWe2001}. Traditionally in the linear space setting, \eqref{ldexact} is approximated with a quadrature rule
	\begin{align*}
    	L_d(q_0,q_1) 
    	&= h\sum\limits_{i=1}^{s} b_i L(q(c_ih),\dot{q}(c_ih))= h \sum\limits_{i=1}^{s} L(Q_i,\dot{Q}_i),
	\end{align*}
and $q(t)$ is approximated by polynomials with degree less than or equal to $s$ with fixed endpoints $q_0$ and $q_1$. By choosing interpolation nodes $\{d_0^\nu\}_{\nu=0}^s$ with $d_0^0=0, d_0^s=1$ and interpolation values $\{q_0^\nu\}_{\nu=0}^s$ with $q_0^0 = q_0$ and $q_0^s=q_1$, $q(t)$ can be expressed as $q(t) = \sum\nolimits_{\nu=0}^s q_0^ \nu \phi_ \nu \left(\frac{t}{h} \right)$ on $[0,h]$, where $\phi_\nu(t)$ are Lagrange polynomials corresponding to the nodes $\{d_0^\nu\}_{\nu=0}^s$. By taking variations with respect to the interpolation values $\{q_0^\nu\}_{\nu=1}^{s-1}$, $q(t)$ is varied over the finite-dimensional function space,
	\[
		M^s = \{q(t)\mid q(t)\in P_s[0,h], q(0)=q_0,q(h)=q_1\}.
	\]
	Consider the quadrature approximation of the action integral, viewing it as a function of the endpoint and interpolation values,
	\[
		F(q_0,q_1,\{q_0^\nu\}_{\nu=1}^{s-1}) = h\sum\limits_{i=1}^{s}b_iL(q(c_ih),\dot{q}(c_ih)),
	\]
	where $q(t) = \sum\nolimits_{\nu=0}^sq_0^\nu\phi_\nu(\frac{t}{h})$. Then, a variational integrator \eqref{vi} can be obtained as follows,
	\begin{equation}\label{gvi}
		\left\{
		\begin{aligned}
			0 &= \frac{\partial F}{\partial q_0^\nu},\qquad \nu=1,2\dots s-1, \\
			-p_0 &= \frac{\partial F}{\partial q_0}, \\
			p_1 &= \frac{\partial F}{\partial q_1}.
		\end{aligned}
		\right.
	\end{equation}
	However, \eqref{gvi} is often impractical due to the complexity of evaluating $q(c_ih)$ and $\dot{q}(c_ih)$, which involve the Lagrange interpolating polynomials. In addition, computing the root of a system of nonlinear equations in \eqref{gvi} can be challenging because the formulation of a fixed-point problem could be complicated, and convergence issues could arise. In contrast, Runge-Kutta methods are already in fixed-point formulation and are convergent as long as consistency is satisfied. 
	
	Now, note that the finite-dimensional function space $M^s$ does not depend on the choice of nodes $\{d_0^\nu\}_{\nu=1}^{s-1}$, and by a tricky elimination of $\phi_\nu(t)$, \eqref{gvi} can be simplified to yield
	\begin{subequations}\label{sprkL}
	\begin{align}
		q_1 &= q_0 + h\sum\limits_{i=1}^sb_i\dot{Q}_i,
			&\quad	p_1&= p_0+h\sum\limits_{i=1}^sb_i\dot{P}_i,\\
		Q_i &= q_0 + h\sum\limits_{j=1}^sa_{ij}\dot{Q}_j,
			&\quad	P_i &= p_0 + h\sum\limits_{j=1}^s\tilde{a}_{ij}\dot{P}_j,\\
		P_i &= \frac{\partial L}{\partial \dot{q}}(Q_i,\dot{Q}_i),
			&\quad	\dot{P}_i &= \frac{\partial L}{\partial q}(Q_i,\dot{Q}_i),
	\end{align}
	\end{subequations}
	where $\tilde{a}_{ij} = b_j(1-\frac{a_{ji}}{b_i})$. When transformed to the Hamiltonian side, \eqref{sprkL} simply recovers the symplectic partitioned Runge--Kutta method. 
	
	The same variational integrator can be derived in a much simpler way: Instead of performing variations on internal points $\{q_0^\nu\}_{\nu=1}^{s-1}$, we will perform variations on the internal derivatives $\{\dot{Q}\}_{i=1}^s$, subject to the constraint $q_1 = q_0 + h\sum\nolimits_{i=1}^s b_i \dot{Q}_i$. Then, the internal points  are reconstructed using the fundamental theorem of calculus and the degree of precision of the quadrature rule to obtain $Q_i = q_0 + h \sum\nolimits_{j=1}^s a_{ij} \dot{Q}_j$. Now, consider the quadrature approximation of the action integral, viewed as a function of the endpoint values and the internal velocities,
	\[
		\tilde{F}(q_0, q_1, \{\dot{Q}_i\}_{i=1}^s, \lambda) = h \sum\limits_{i=1}^s b_i L(Q_i,\dot{Q}_i) + \lambda^{T} \left(q_1 - q_0 - h \sum\nolimits_{i=1}^s b_i\dot{Q}_i\right),
	\]
	where $\lambda$ is a Lagrange multiplier that enforces the constraint. Then, a variational integrator \eqref{vi} can be obtained as follows,
	\begin{equation}
	\label{ngvi}
		\left\{
		\begin{aligned}
			0 &= \frac{\partial \tilde{F}}{\partial \dot{Q}_i},\qquad i = 1, 2 \ldots s, \\
			0 &= \frac{\partial \tilde{F}}{\partial \lambda}, \\
			Q_i &= q_0 + h\sum\limits_{j=1}^sa_{ij}\dot{Q}_j,\\
			-p_0 &= \frac{\partial \tilde{F}}{\partial q_0}, \\
			p_1 &= \frac{\partial \tilde{F}}{\partial q_1}.
		\end{aligned}
		\right.
	\end{equation}
	Explicitly expanding \eqref{ngvi} and eliminating the Lagrange multiplier yields \eqref{sprkL} in a much more straightforward manner. This same technique, known as the \textit{constrained Galerkin method}, is adopted on the rotation group $SO(n)$ in order to directly obtain a variational integrator in fixed-point form.

\section{Euler--Poincar\'e \& Lie--Poisson Reductions}
\label{append_LPReductions}
	When the Lagrangian $L$ or Hamiltonian $H$ is left-invariant, the mechanical system can be symmetry reduced to evolve on the Lie algebra $\mathfrak{g}$ or its dual $\mathfrak{g}^*$, respectively, assuming some regularity. On the Lagrangian side, this corresponds to Euler--Poincar\'e reduction~\cite{MaRa1999}, which is described by the Euler--Poincar\'e equations,
	\[
		\frac{d}{dt}\left(\frac{\partial \bm{l}}{\partial \epsilon}\right) 
		= \ad\nolimits_{\epsilon}^*\left(\frac{\partial \bm{l}}{\partial \epsilon}\right).
	\]
	The above is expressed in terms of the reduced Lagrangian $\bm{l}(\epsilon) = \bm{l}(g^{-1}\dot{g})=L(g,\dot g)$. As a result, this can be described in terms of a reduced variational principle $\delta \int_a^b \bm{l}(\epsilon(t)) \, dt = 0$ with respect to constrained variations of form $\delta \epsilon = \dot{\eta}+[\epsilon,\eta]$, where $\eta(t)$ is an arbitrary path in the Lie algebra $\mathfrak{g}$ that vanishes at the endpoints, namely $\eta(a)=\eta(b)=0$. The constraint on the variations $\delta \epsilon$ are induced by the condition that $\epsilon=g^{-1}\dot{g}$ together with arbitrary variations $\delta g$ that vanish at the endpoints.
	
	On the Hamiltonian side, this corresponds to Lie--Poisson reduction~\cite{MaRa1999}. Recall that the Lie--Poisson structure on $\mathfrak{g}^*$ is given by
	\[
		\{F,G\}(\mu) = \left\langle\mu, \left[\frac{\partial F}{\partial \mu}, \frac{\partial G}{\partial \mu}\right]\right\rangle,
		\]
	and together with the reduced Hamiltonian $\bm{h}(\mu)$, they gives the Lie--Poisson equations on $\mathfrak{g}^*$,
	\[
		\frac{d\mu}{dt} = \text{ad}_{\frac{\partial \bm{h}}{\partial \mu}}^*(\mu).
	\]
	If the discrete Lagrangian $L_d(g_0,g_1)$ is also $G$-invariant, meaning $L_d(g\cdot g_0,g\cdot g_1) = L_d(g_0,g_1)$ for some $g \in G$, then \eqref{vi} can be reduced to a Lie--Poisson integrator~\cite{Pekarsky1999},
	\begin{equation}
	\label{liepoisson}
		\left\{
		\begin{aligned}
			{\mu}_0 &= \bm{l}_d^{'} (f_0)f_0^{-1}, \\
			{\mu}_1 &= f_0^{-1}\cdot {\mu}_0\cdot f_0,
		\end{aligned}
		\right.
	\end{equation}
where $\bm{l}_d(f_0) = L_d(e,f_0)$. This algorithm preserves the coadjoint orbits and, hence, the Poisson structure on $\mathfrak{g}^*$.

\section{Detailed Derivations for the Lagrangian Variational Integrators}
\label{append_detailedComp}
	\subsection{Derivations of $\frac{\partial \tilde{F}}{\partial \Omega_k}$}
	\label{append_dFdOmega}
	Initially, we have 
	\begin{align*}
		\left. \frac{d} {d t} \right|_{t=0} \tilde{F}(g_0, g_1, \ldots, \Omega_k(t),  \ldots, \Lambda) 
		&= h \sum\limits_{i=1}^s b_i \left\langle \frac{\partial L}{\partial U} (U_i,\Omega_i), X_{ik} \right\rangle + h b_k \left\langle \frac{\partial L}{\partial \Omega} (U_k,\Omega_k), \delta\Omega_k \right\rangle \\
    		&	\qquad + \left\langle \Lambda, \asym \left(h g_1^{T} \sum\nolimits_{i=1}^s b_i U_i X_{ik} \Omega_i + h g_1^{T} b_k U_k \delta\Omega_k \right)\right\rangle,
	\end{align*}
where we can use the properties of the inner products to express the last term as follows,
	\[
		 h \sum\limits_{i=1}^s b_i \langle \asym (U_i^{T} g_1 \Lambda \Omega_i^{T}), X_{ik} \rangle + h b_k \langle \asym(U_k^{T} g_1 \Lambda),\delta\Omega_k \rangle.
	\]
	Then, we continue using equation \eqref{Xik} to obtain,
	\begin{align*}
		\left. \frac{d} {d t} \right|_{t=0} \tilde{F}(g_0, g_1, \ldots, \Omega_k(t),  \ldots, \Lambda)
		&= h \sum\limits_{i=1}^s b_i \left\langle \frac{\partial L}{\partial U} (U_i, \Omega_i) + \asym (U_i^{T} g_1 \Lambda\Omega_i^{T}), X_{ik} \right\rangle \\
			& \qquad + h b_k \left\langle \frac{\partial L}{\partial \Omega}(U_k,\Omega_k) + \asym (U_k^{T} g_1 \Lambda), \delta\Omega_k \right\rangle \\
		&= h \sum\limits_{i=1}^s b_i \left\langle \frac{\partial L}{\partial U} (U_i, \Omega_i) + \asym (U_i^{T} g_1 \Lambda\Omega_i^{T}), h (\tilde{\pi}_i \circ \tilde{\phi}^{-1}\circ \tilde{\psi}_k)(\delta\Omega_k) \right\rangle \\
			& \qquad + h b_k \left\langle \frac{\partial L}{\partial \Omega}(U_k,\Omega_k) + \asym (U_k^{T} g_1 \Lambda), \delta\Omega_k \right\rangle \\
		&= h^2 \sum\limits_{i=1}^s b_i \left\langle (\tilde{\pi}_i \circ \tilde{\phi}^{-1}\circ \tilde{\psi}_k)^* \left( \frac{\partial L}{\partial U} (U_i, \Omega_i) + \asym (U_i^{T} g_1 \Lambda\Omega_i^{T}) \right), \delta\Omega_k \right\rangle \\
			& \qquad + h b_k \left\langle \frac{\partial L}{\partial \Omega}(U_k,\Omega_k) + \asym (U_k^{T} g_1 \Lambda), \delta\Omega_k \right\rangle.
	\end{align*}
	We finally have
	\[
		\frac{\partial \tilde{F}}{\partial \Omega_k}
		= 
		h^2 \sum\limits_{i=1}^s b_i (\tilde{\pi}_i \circ \tilde{\phi}^{-1}\circ \tilde{\psi}_k)^* \left( \frac{\partial L}{\partial U} (U_i, \Omega_i) + \asym (U_i^{T} g_1 \Lambda\Omega_i^{T}) \right) +  h b_k \left( \frac{\partial L}{\partial \Omega}(U_k,\Omega_k) + \asym (U_k^{T} g_1 \Lambda) \right).
	\]

	\subsection{Explicit Expression for $\tilde{\phi}^*$: A Derivation}
	\label{append_phistar}
	Let us consider $(S_1, \ldots, S_s), (\tilde{S}_1, \ldots, \tilde{S}_s) \in \sk (n)^s$, and so 
	\begin{align*}
		\langle \tilde{\phi}^*(S_1, \ldots, S_s), (\tilde{S}_1, \ldots, \tilde{S}_s) \rangle
		&= \langle (S_1, \ldots, S_s), \tilde{\phi}(\tilde{S}_1, \ldots, \tilde{S}_s) \rangle	\\
		&= \left\langle (S_1, \ldots, S_s), \left\{ \tilde{S}_i - d\mathbb{P}_{A_i} \left( h \sum\nolimits_{j=1}^s a_{ij} U_j \tilde{S}_j \Omega_j \right) \right\}_{i=1}^s \right\rangle	\\
		&= \sum\limits_{i=1}^s \left\langle S_i, \tilde{S}_i - d\mathbb{P}_{A_i} \left( h \sum\nolimits_{j=1}^s a_{ij} U_j \tilde{S}_j \Omega_j \right) \right\rangle \\
		&= \sum\limits_{i=1}^s \langle S_i,\tilde{S}_i \rangle - \sum\limits_{i=1}^s \tr \left(d \mathbb{P}_{A_i}^*(S_i) \left(h \sum\nolimits_{j=1}^s a_{ij} U_j \tilde{S}_j \Omega_j\right)^T \right) \\
		&= \sum\limits_{i=1}^s \langle S_i,\tilde{S}_i \rangle - \sum\limits_{i,j=1}^s \tr \left( \left( h a_{ij} U_j^T d\mathbb{P}_{A_i}^*(S_i) \Omega_j^T\right) \tilde{S}_j^T \right) \\
		&= \sum\limits_{i=1}^s \langle S_i,\tilde{S}_i \rangle - \sum\limits_{j=1}^s \tr \left(\left( h U_j^T \sum\nolimits_{i=1}^s a_{ij} d\mathbb{P}_{A_i}^*(S_i) \Omega_j^T\right) \tilde{S}_j^T \right)	\\
		&= \sum\limits_{j = 1} \left \langle S_j - \asym \left( h U_j^T \sum\nolimits_{i=1}^s a_{ij} d\mathbb{P}_{A_i}^*(S_i) \Omega_j^T \right) , \tilde{S}_j \right\rangle.
	\end{align*}
	This gives us equation \eqref{phistar}.
	
	\subsection{Explicit Expression for $\tilde{\psi}_k^*$: A Derivation}
	\label{append_psistar}
	Consider $(S_1, \ldots, S_s),\in \sk (n)^s$ and $\tilde{S} \in \sk (n)$, then
	\begin{align*}
		\langle \tilde{\psi}_k^* (S_1, \ldots, S_s), \tilde{S} \rangle
		&= \langle (S_1, \ldots, S_s), \tilde{\psi}_k(\tilde{S}) \rangle	\\
		&= \left\langle (S_1, \ldots, S_s), \left\{ d \mathbb{P}_{A_i} (a_{ik} U_k \tilde{S})\right\}_{i=1}^s \right\rangle \\
    	&= \sum\limits_{i=1}^s  \tr \left( d \mathbb{P}_{A_i}^*(S_i) (a_{ik} U_k \tilde{S})^T \right) \\
    	&= \sum\limits_{i=1}^s  \tr \left( (a_{ik} U_k^T d \mathbb{P}_{A_i}^*(S_i)) \tilde{S}^T \right)	\\
    	&= \left \langle \asym \left( U_k^T \sum\nolimits_{i=1}^s a_{ik} d \mathbb{P}_{A_i}^*(S_i) \right), \tilde{S} \right \rangle.
	\end{align*}
	This gives us equation \eqref{psistar}.
	
	\subsection{Explicit Expression for $\tilde{\varphi}^*$: A Derivation}
	\label{append_varphistar}
	Let us derive $\tilde{\varphi}^*$ by considering $(S_1, \ldots, S_s)\in \sk (n)^s$ and $\tilde{S} \in \sk (n)$, and so
	\begin{align*}
		\langle\tilde{\varphi}^*(S_1, \ldots, S_s), \tilde{S} \rangle 
		&= \left\langle(S_1, \ldots, S_s), \{d \mathbb{P}_{A_i}(g_0 \tilde{S})\}_{i=1}^s \right\rangle \\
		&= \sum\limits_{i=1}^s \tr \left(d \mathbb{P}_{A_i}^*(S_i)(g_0 \tilde{S})^T \right) \\
		&= \sum\limits_{i=1}^s \tr \left( (g_0^{T} d \mathbb{P}_{A_i}^*(S_i)) \tilde{S}^T \right) \\
		&= \left\langle \asym \left(g_0^{T}\sum\nolimits_{i=1}^s d \mathbb{P}_{A_i}^*(S_i) \right), \tilde{S} \right\rangle.
	\end{align*}
	This gives us equation \eqref{varphistar}.
	
	\subsection{Derivation of $\frac{\partial \tilde{F}}{\partial g_0}$}
	\label{append_dFdg0}
	We compute
	\begin{align*}
    	\left.\frac{d}{d t}\right|_{t=0} \tilde{F}(g_0(t), g_1,\{ \Omega_i \}_{i=1}^s ,\Lambda)
		&= h \sum\limits_{i=1}^s b_i \left\langle \frac{\partial L}{\partial U}(U_i,\Omega_i), Y_i \right\rangle + \left\langle \Lambda, \asym \left( g_1^{T}\left( g_0 \delta g_0 + h \sum\nolimits_{i=1}^s b_i U_i Y_i\Omega_i \right)\right) \right\rangle	\\
		&= h \sum\limits_{i=1}^s b_i \left\langle \frac{\partial L}{\partial U}(U_i,\Omega_i), Y_i \right\rangle + \tr \left( \Lambda (g_1^T g_0 \delta g_0)^T \right) + h \sum\limits_{i=1}^s \tr \left( \Lambda (b_i g_1^T U_i Y_i \Omega_i)^T \right)	\\
		&= h \sum\limits_{i=1}^s b_i \left\langle \frac{\partial L}{\partial U}(U_i,\Omega_i), Y_i \right\rangle + \tr \left( (g_0^T g_1\Lambda )\delta g_0^T \right) + h \sum\limits_{i=1}^s \tr \left(b_i(U_i^T g_1 \Lambda \Omega_i^T)Y_i \right)	\\
		&= h \sum\limits_{i=1}^s b_i \left\langle \frac{\partial L}{\partial U}(U_i,\Omega_i) + \asym (U_i^{T} g_1 \Lambda \Omega_i^{T}), Y_i \right\rangle + \langle\asym (g_0^{T} g_1 \Lambda), \delta g_0 \rangle	\\
		&= h \sum\limits_{i=1}^s b_i \left\langle \frac{\partial L}{\partial U}(U_i,\Omega_i) + \asym (U_i^{T} g_1 \Lambda \Omega_i^{T}), (\tilde{\pi}_i \circ \tilde{\phi}^{-1} \circ \tilde{\varphi}) (\delta g_0) \right\rangle 	\\
			& \qquad + \langle\asym (g_0^{T} g_1 \Lambda), \delta g_0 \rangle	\\
	\end{align*}
	Thus, we have 
	\[
		\frac{\partial \tilde{F}}{\partial g_0}
		=  h \sum\limits_{i=1}^s b_i (\tilde{\pi}_i \circ \tilde{\phi}^{-1} \circ \tilde{\varphi})^* \left( \frac{\partial L}{\partial U}(U_i,\Omega_i) + \asym (U_i^{T} g_1 \Lambda \Omega_i^{T}) \right) + \asym (g_0^{T} g_1 \Lambda).
	\]

\bibliography{polardecompostion}
\bibliographystyle{plainnat}

\end{document}